\documentstyle[amscd,amssymb,verbatim,12pt]{amsart}
\input{diagrams}

\newcommand{\tr}{\operatorname{tr}}

\renewcommand{\mod}{\operatorname{mod}}

\newcommand{\Cone}{\operatorname{Cone}}

\newcommand{\Tors}{\operatorname{Tors}}

\newcommand{\OO}{{\cal O}}

\newcommand{\Vect}{\operatorname{Vect}}

\newcommand{\coker}{\operatorname{coker}}

\newcommand{\MM}{{\cal M}}

\newcommand{\BB}{{\cal B}}

\newcommand{\SL}{\operatorname{SL}}

\newcommand{\hra}{\hookrightarrow}
\newcommand{\lan}{\langle}
\newcommand{\ran}{\rangle}
\newcommand{\Coh}{\operatorname{Coh}}

\newcommand{\CC}{{\cal C}}

\renewcommand{\P}{{\Bbb P}}

\newcommand{\Ind}{\operatorname{Ind}}
\newcommand{\ga}{\gamma}
\newcommand{\de}{\delta}
\newcommand{\eps}{\epsilon}

\renewcommand{\ker}{\operatorname{ker}}
\newcommand{\im}{\operatorname{im}}

\numberwithin{equation}{subsection}

\newtheorem{thm}{Theorem}[subsection]
\newtheorem{prop}[thm]{Proposition}
\newtheorem{prop-def}[thm]{Proposition-Definition}
\newtheorem{lem}[thm]{Lemma}
\newtheorem{cor}[thm]{Corollary}
\newenvironment{rem}{\vspace{3mm}\noindent
{\bf Remark.}}{\vspace{3mm}}
\newenvironment{defi}{\vspace{3mm}\noindent
{\bf Definition.}}{\vspace{3mm}}
\newenvironment{rems}{\vspace{3mm}
\noindent {\bf Remarks.}}{\vspace{3mm}}
\newenvironment{ex}{\vspace{3mm}\noindent
{\bf Example.}}{\vspace{3mm}}

\newcommand{\Pf}{\noindent {\it Proof}}

\newcommand{\ov}{\overline}

\renewcommand{\Im}{\operatorname{Im}}

\newcommand{\rk}{\operatorname{rk}}

\renewcommand{\AA}{{\cal A}}
\newcommand{\FF}{{\cal F}}
\newcommand{\TT}{{\cal T}}

\newcommand{\PP}{{\cal P}}

\newcommand{\Hom}{\operatorname{Hom}}
\newcommand{\Ext}{\operatorname{Ext}}
\newcommand{\End}{\operatorname{End}}

\renewcommand{\a}{\alpha}
\renewcommand{\b}{\beta}

\newcommand{\th}{\theta}
\newcommand{\C}{{\Bbb C}}
\newcommand{\R}{{\Bbb R}}
\newcommand{\Z}{{\Bbb Z}}

\newcommand{\Ga}{\Gamma}

\newcommand{\wt}{\widetilde}

\newcommand{\sign}{\operatorname{sign}}

\newcommand{\sub}{\subset}
\newcommand{\ed}{\qed\vspace{3mm}}

\newcommand{\Qcoh}{\operatorname{Qcoh}}

\newcommand{\dlim}{{\displaystyle \lim_{\longrightarrow}}}

\title{Quasicoherent sheaves on complex noncommutative two-tori}
\author{A. Polishchuk}
\thanks{Supported in part by NSF grant}

\begin{document}
\begin{abstract}  We introduce the notion of a {\it quasicoherent sheaf} on a complex noncommutative two-torus $T$ as an ind-object in the category of holomorphic vector bundles on $T$. 
Extending the results of \cite{PS} and \cite{P-holbun}
we prove that the derived category of quasicoherent sheaves on $T$ is equivalent to the derived
category of usual quasicoherent sheaves on the corresponing elliptic curve. We define
the rank of a quasicoherent sheaf on $T$ that can take arbitrary nonnegative
real values. We study the category
$\Qcoh(\eta_T)$ obtained by taking the quotient of the category of quasicoherent sheaves by
the subcategory of objects of rank zero (called {\it torsion sheaves}). We show that
projective objects of finite rank in $\Qcoh(\eta_T)$ are classified up to an isomorphism by their rank.
We also prove
that the subcategory of objects of finite rank in $\Qcoh(\eta_T)$ is equivalent to the category of
finitely presented modules over a semihereditary algebra.
\end{abstract}
\maketitle

\bigskip

\centerline{\sc Introduction}

\medskip

The goal of this paper is to define and study the category of quasicoherent
sheaves on a noncommutative two-torus equipped with a complex structure.
Recall that a noncommutative two-torus $T_{\th}$ is defined via its algebra
of smooth functions $A_{\th}$ that is determined by an irrational real number
$\th$. A complex structure on $T_{\th}$ is given by a certain derivation
$\de_{\tau}$ of the algebra $A_{\th}$ associated with a complex parameter $\tau$
(see section \ref{prelim-sec}). We view this derivation as an analogue of the operator
$\ov{\partial}$. We denote by $T=T_{\th,\tau}$ the
obtained complex noncommutative torus.
Holomorphic vector bundles on $T$
are defined as finitely generated projective
right $A_{\th}$-modules equipped with a lifting of $\de_{\tau}$
(see section \ref{prelim-sec}). 

The category $\Vect(T)$
of holomorphic vector bundles on $T$
was studied in \cite{PS}, \cite{P-real} and \cite{P-holbun}.
In particular, it was proved in \cite{P-holbun} that this category is abelian. Furthermore, 
if one tries to mimick the usual definition of a coherent sheaf in this
situation one obtains that every such coherent sheaf is a vector bundle (see
Theorem \ref{coh-thm}). 
Thus, there are no analogues of coherent torsion sheaves
in our situation. However, as we will show, things become more interesting
if we consider quasicoherent sheaves.
We define the category $\Qcoh(T)$
of quasicoherent sheaves as the category of ind-objects
in $\Vect(T)$. We realize this category explicitly 
as a full subcategory in the
larger category of {\it holomorphic modules} on $T$ (these are arbitrary
$A_{\th}$-modules equipped with a lifting of $\de_{\tau}$,
see section \ref{quasi-adm-sec}).

Our first result is that the category $\Qcoh(T)$ is equivalent to a certain abelian subcategory in
the derived category $D^b(\Qcoh(E))$ of usual quasicoherent sheaves on the elliptic
curve $E=\C/(\Z+\Z\tau)$. This extends a similar result for $\Vect(T)$ (see \cite{PS}, \cite{P-holbun}). Furthermore,
there is an equivalence of derived categories $D^b(\Qcoh(T))\simeq D^b(\Qcoh(E))$ (see
Theorem \ref{FM-tori-thm}).

Next, we prove that
the rank of a vector bundle on $T$ (which
is a nonnegative real number of the form $m\th+n$ with $m,n\in\Z$) extends
naturally to an additive function on the category of quasicoherent sheaves taking
all possible values in $\R_{\ge 0}\cup\{+\infty\}$
(see Theorem \ref{rank-thm} and Corollary \ref{quasi-main-cor}).

We define torsion sheaves as quasicoherent sheaves of rank zero.
It turns out that there are many nontrivial torsion sheaves on $T$.
In fact, the subcategory $\Tors$ of torsion sheaves is big enough 
to make the passage from $\Qcoh(T)$ to
the quotient-category 
$$\Qcoh(\eta_T):=\Qcoh(T)/\Tors$$ 
in many ways similar to the passage to the general
point in commutative algebraic geometry. 

Note that by definition the rank descends to
an additive function on $\Qcoh(\eta_T)$.
Let $\Qcoh^f(\eta_T)$ denote the full subcategory of $\Qcoh(\eta_T)$ consisting
of quasicoherent sheaves of finite rank.

Our main result is that the category $\Qcoh^f(\eta_T)$ is equivalent to the category of finitely presented modules over a semihereditary ring and that projective objects of $\Qcoh^f(\eta_T)$ are uniquely
determined by their rank that can be an arbitrary real number (see Theorem \ref{cat-eq-thm},
Proposition \ref{semiher-prop},
Corollaries \ref{proj-isom-cor} and \ref{K0-cor}).
We also show that projective objects of $\Qcoh^f(\eta_T)$ correspond to {\it quasi vector
bundles}, i.e., filtering unions of holomorphic vector bundles (see Theorem \ref{proj-thm}).
The above equivalence of categories is obtained by taking 
a projective object $P\in\Qcoh^f(\eta_T)$ and sending $X\in\Qcoh^f(\eta_T)$ to
the module $\Hom_{\Qcoh(\eta_T)}(P,X)$ over the ring 
$R_P=\End_{\Qcoh}(\eta_T)(P)$.

To a large extent the study of the category $\Qcoh(\eta_T)$ reduces to the problem of
constructing holomorphic subbundles and quasicoherent subsheaves
in a stable holomorphic bundle $V$ on $T$ with given properties (see section \ref{prelim-sec}
for the definition of stability for bundles on $T$). 
One of the constructions we use can be considered as a categorification of a 
two-sided version of the continuous fraction process (see sections
\ref{div-constr-sec} and \ref{subbun-sec}). Subbundles constructed in
this way are numbered by vertices of a binary tree and depend also on
some continuous parameters. Theorem \ref{main-thm} implies that
the ranks of the subbundles of $V$ obtained by this construction
are constrained only by the requirement that the corresponding
slope is smaller than the slope of $V$ (of course, these ranks also
should be smaller than $\rk V$).

One may wonder what kind of restrictions one gets for $A_{\th}$-modules underlying
quasicoherent sheaves on $T$. 
We consider the subcategory of {\it countably generated} quasicoherent
sheaves and show that it is a Serre subcategory in $\Qcoh(T)$ and that the underlying $A_{\th}$-modules always have projective dimension $\le 1$
but are not necessarily projective (see Theorem \ref{count-thm}, Corollary \ref{quasi-main-cor}). 
In particular, we derive
that every quasicoherent ideal in $A_{\th}$ is countably generated.
The proof is based on the analogue of the Harder-Narasimhan filtration for quasicoherent subsheaves
of holomorphic vector bundles on $T$ (see section \ref{HN-sec}).

One may view the ring $R_P$ appearing above as an algebraic version of
the von Neumann factor of type $II_1$. Namely, if $P$ is a vector bundle then
such a factor appears as the closure of the endomorphism
algebra $\End_{A_{\th}}(P)$ in the algebra of bounded operators on the appropriate
Hilbert space. We conjecture that algebra $R_P$ contains a ``convergent"
subalgebra $R'_P$ with $K_0(R'_P)=K_0(R_P)=\R$ such that $R'_P$ embeds
also into the above von Neumann factor.

The plan of the paper is as follows. In section \ref{facts-sec} we prove some auxiliary statements
about the category $\Vect(T)$ of holomorphic vector bundles. 
In section \ref{quasi-sec} we define and study the category $\Qcoh(T)$ of quasicoherent sheaves.
This includes realizing $\Qcoh(T)$ explicitly as the category of admissible holomorphic modules
in \ref{quasi-adm-sec}, establishing the equivalence of derived categories in \ref{FM-sec}, and constructing the extension of the rank function to quasicoherent sheaves in \ref{rank-sec}. We also study 
the subcategory of countably presented quasicoherent sheaves in \ref{count-sec}.
Finally, in section \ref{quasi-gen-sec} we prove our results about the category $\Qcoh(\eta_T)$ of 
``sheaves at the general point of $T$" (and its subcategory $\Qcoh^f(\eta_T)$).

\noindent
{\it Acknowledgment}. I am grateful to Paul Smith for the stimulating question on possible ranks
of holomorphic ideals in $A_{\th}$. 
The answer to this question is Theorem \ref{main-thm}.
I'm also indebted to Amnon Neeman for a helpful discussion on inductive
limits in derived categories. 

\section{Some facts about the category of holomorphic bundles on $T$}\label{facts-sec}

In this section we prove some results about holomorphic bundles on a
noncommutative torus $T$ that will be used in our study of quasicoherent sheaves on $T$.
After providing some background we describe in \ref{subbun-sec} two constructions of subbundles
in holomorphic vector bundles on $T$ that will play a crucial role in section \ref{quasi-gen-sec}.
Then in \ref{deform-sec} we show that every two holomorphic vector bundles of the same rank on $T$ are deformation equivalent.

\subsection{Preliminaries}\label{prelim-sec}

Throughout this paper the number $\th$ is assumed to be irrational. 
By a module over a ring we always mean a right module (same convention
for ideals).

The algebra $A_{\th}$ of smooth functions on the noncommutative torus $T_{\th}$
is defined as the algebra of series of the form 
$\sum_{(m,n)\in\Z^2}a_{m,n}U_1^mU_2^n$,
where the generators $U_1$ and $U_2$ satisfy the commutation relation 
$U_1U_2=\exp(2\pi i \th)U_2U_1$, and 
$(a_{m,n})$ is a collection of complex numbers rapidly decreasing as $m^2+n^2\to\infty$.
 We fix $\tau\in\C$ such that $\Im\tau<0$ and consider the derivation
$$\de_{\tau}:A_{\th}\to A_{\th}: \sum a_{m,n} U_1^mU_2^n\mapsto 2\pi i\sum a_{m,n} (m\tau+n)U_1^mU_2^n$$
as an analogue of the $\ov{\partial}$-operator 
giving the complex structure on our noncommutative torus.
By definition, a {\it vector bundle} on $T_{\th}$ is a finitely generated projective right 
$A_{\th}$-module.
A {\it holomorphic vector bundle} on $T_{\th}$ is a vector bundle $P$ equipped with
an operator $\ov{\nabla}:P\to P$ such that
\begin{equation}\label{Leibnitz-eq}
\ov{\nabla}(sa)=\ov{\nabla}(s)a+s\de_{\tau}(a)
\end{equation}
for all $s\in P$, $a\in A_{\th}$.

The category $\Vect(T)$ of holomorphic vector bundles on a noncommutative
complex torus $T=T_{\th,\tau}$ was studied in \cite{PS} and \cite{P-holbun}. 
We showed that there is an equivalence of $\Vect(T)$
with a certain abelian subcategory $\CC^{\th}$
in the derived category $D^b(E)$
of coherent sheaves on the elliptic curve $E=\C/(\Z+\Z\tau)$.
Nonzero objects $F\in\CC^{\th}$ satisfy 
$$\rk_{\th}(F):=\deg_E(F)\th+\rk_E(F)>0,$$
where $\deg_E$ and $\rk_E$ are the standard degree and rank functions
on $D^b(E)$. On the other hand, for a vector bundle on a noncommutative
torus there is also a notion of rank defined using the trace
$\tr:A_{\th}\to\R:\sum_{m,n}a_{m,n}U_1^mU_2^n\mapsto a_{0,0}$.
For $V\in\Vect(T)$ we denote this rank by $\rk V$. If we view $V$ as an object of $\CC^{\th}$
via the equivalence $\Vect(T)\simeq\CC^{\th}$ then we have $\rk(V)=\rk_{\th}(V)$.
We will also use the functions $\deg(V)=\deg_E(V)$ (degree) and
$\mu(V)=\deg(V)/rk(V)$ (slope). Note that the degree is determined by the rank:
$\deg(V)=m$, where $\rk(V)=m\th+n$ with $m,n\in\Z$.

\begin{defi} We say that an object
$A$ in $D^b(E)$ is {\it stable}
(resp., {\it semistable}) if $A\simeq \FF[n]$, where $\FF$ is either a stable
vector bundle on $E$ or the structure sheaf of a point in $E$ (resp., $\FF$ is
either a semistable vector bundle or a torsion sheaf).
\end{defi}

Note that every object of $D^b(E)$ can be decomposed into a direct sum of semistable
objects. Viewing $\Vect(T)$ as a subcategory in $D^b(E)$ we obtain the definition
of stability and semistability for holomorphic vector bundles on $D^b(E)$.
It is easy to see that stable holomorphic bundles on $T$ correspond exactly to standard
holomorphic structures on basic projective modules over $A_{\th}$ (see \cite{PS}). 
In the following lemma we check that the above definition coincides with the notion of stability
obtained using slopes of bundles on $T$.

\begin{lem}\label{stable-bun-lem} 
A holomorphic vector bundle $V$ on $T$ is stable (resp.,
semistable) iff for every subbundle $W\sub V$ such that $0<\rk W<\rk V$,
one has $\mu(W)<\mu(V)$ (resp., $\mu(W)\le\mu(V)$).
\end{lem}

\Pf . We will only prove the part concerning semistability and leave the stability part to the reader.
Assume first that $V$ is semistable. If $W$ is a semistable bundle
on $T$ then $\Hom(W,V)\neq 0$ only if $\mu(W)\le\mu(V)$ (this follows from Lemma 1.6 of \cite{P-real}
using Serre duality). Since every vector bundle of slope $\mu$ contains a semistable subbundle of slope 
$\ge\mu$, the ``only if" part follows.
Conversely, assume that for every $W\sub V$ with $0<\rk W<\rk V$
one has $\mu(W)\le\mu(V)$. Let $V_0\sub V$ be a maximal semistable subbundle of maximal slope (it
exists). Then $\mu(V_0)\le\mu(V)$ hence we should have $V_0=V$, i.e., $V$ is semistable.
\ed

\begin{lem}\label{triangle-lem} 
Let $X\to Y\to Z\to X[1]$ be an exact triangle in $D^b(E)$ with $X,Y\in\CC^{\th}$.
Assume that $Z=\oplus_{i=1}^k Z_i$, where each $Z_i$ is semistable with
$\rk_{\th}(Z_i)>0$ (resp., $\rk_{\th}(Z_i)<0$). Then $Z\in\CC^{\th}$ (resp., $Z[-1]\in\CC^{\th}$).
\end{lem}

\Pf . Let $H^i:D^b(E)\to\CC^{\th}$ be the cohomology functors associated with
the $t$-structure that has $\CC^{\th}$ as a heart. Then we have
$Z\simeq H^0(Z)\oplus H^{-1}(Z)[1]$. Now the condition that $\rk_{\th}$ takes
positive (resp., negative) values on semistable summands
of $Z$ implies that $H^{-1}(Z)=0$ (resp., $H^0(Z)$).
\ed

\begin{rem}
Statements similar to Lemmas \ref{stable-bun-lem} and \ref{triangle-lem} 
hold in the more general framework of stability 
conditions on derived categories developed in \cite{Bridge}. 
More precisely, in 
Lemma \ref{triangle-lem} one has to replace the rank with
the imaginary part of the central charge function. Also, $Z_i$'s 
should be replaced by the Harder-Narasimhan constituents of $Z$. 
\end{rem}

Recall that if $\th$ and $\th'$ are related by a fractional-linear transformation
then the algebras $A_{\th}$ and $A_{\th'}$ are
Morita equivalent. The corresponding categories of holomorphic bundles
$\Vect(T_{\th,\tau})$ and $\Vect(T_{\th',\tau})$ are equivalent. For every such
Morita equivalence $\Phi:\Vect(T)\wt{\to}\Vect(T')$ we have
$\rk(\Phi(V))=c\cdot\rk(V)$ for some constant $c>0$.
Moreover, for every stable vector bundle $V$ on $T$ there exists a
Morita equivalence $\Phi$ such that $\rk(\Phi(V))=1$.

\subsection{Binary division process associated with an irrational number}
\label{div-constr-sec}

Let us denote $L_{\th}=\Z\th+\Z$ and let $\PP\sub L_{\th}$ be the set
of primitive vectors in $L_{\th}$. We also denote $\PP_{>0}=\PP\cap(0,+\infty)$.

We equip $L_{\th}=\Z\th+\Z$ with a $\Z$-valued bilinear form $\chi=\chi_{\th}$ by 
setting
$$\chi(m\th+n,m'\th+n')=m'n-mn'.$$
Note that $L_{-\th}=L_{\th}$ but $\chi_{-\th}=-\chi_{\th}$.
Recall that $\PP\sub L_{\th}$ denotes the set of primitive vectors and
$\PP_{>0}=\PP\cap(0,+\infty)$.

\begin{lem} For every $v\in\PP_{>0}$ there exists a unique vector 
$\phi(v)\in\PP_{>0}$
such that $\phi(v)<v$ and $\chi(\phi(v),v)=1$.
\end{lem}

\Pf . Let $v=m\th+n$. We are looking for a vector $m_1\th+n_1\in L_{\th}$ 
such that $0<m_1\th+n_1<m\th+n$ and
$mn_1-m_1n=1$. Thus, it suffices to prove the existence and uniqueness
of $m_1\in\Z$ such that $m_1n\equiv -1(m)$ and
$$0<\frac{m_1(m\th+n)+1}{m}<m\th+n.$$
The latter condition is equivalent to
\begin{equation}\label{main-ineq}
0<\sign(m)[m_1+(m\th+n)^{-1}]<|m|.
\end{equation}
Thus, $m_1$ should be within a given interval of length $|m|$
with irrational ends. Since the residue of $m_1$ modulo $|m|$ is
fixed by the condition $m_1n\equiv -1(m)$ we get a unique solution.
\ed

Using the map $\phi=\phi_{\th}:\PP_{>0}\to \PP_{>0}$ we can define a canonical way 
to divide every segment
$[a,b]$ such that $b-a\in\PP_{>0}$ into two subsegments: $[a,a+\phi(b-a)]$ and 
$[a+\phi(b-a),b]$.
Moreover, each of the two new segments also has the length in $\PP_{>0}$.
Let us start with the segment $[0,1]$ and divide it into two subsegments using this 
recipe.
Then divide each of the new segments in two subsegments again, etc.
Below we will refer to the subsegments $[a,b]\sub [0,1]$ that are being divided
in this process as division subsegments.
Let $\BB_{\th}\sub (0,1)$ denote the set of endpoints of all the division subsegments.  

\begin{thm}\label{main-div-thm} 
The set $\BB_{\th}$ coincides with the set of all $m\th+n\in L_{\th}\cap (0,1)$ such 
that $m<0$.
\end{thm}

\Pf . It is easy to see that $\chi_{-\th}=\chi_{\th}$ on $L_{-\th}=L_{\th}$. 
Hence, 
$\phi_{-\th}(v)=v-\phi_{\th}(v)$ for $v\in\PP_{>0}$. This implies that $\BB_{-
\th}=1-\BB_{\th}$.
Thus, the assertions of the theorem for $\th$ and for $-\th$ are equivalent.
Let us first prove that for $m\th+n\in\BB_{\th}$ one has $m<0$.
Assume that $v=m\th+n$ is an element of the $k$-th generation
of points obtained by our division process. We use induction in $k$. The first 
point of the division process
is $a-\th$, where $a$ is the unique integer such that $0<a-\th<1$, 
so the assertion holds for $k=1$. Assume that our claim is true for all $k'<k$ and 
for all $\th$.
Changing $\th$ to $-\th$ if necessary we can assume that $v<a-\th$.
Set $\th'=1/(a-\th)$ and $v'=v/(a-\th)$. We claim that $v'$ is a 
$(k-1)$-th generation point of the division
process associated with $\th'$. Indeed, the map $\a:w\mapsto w/(a-\th)$ is an 
order-preserving
isomorphism from $L_{\th}$ to $L_{\th'}$. Furthermore, since $\a(1)=\th'$ and 
$\a(\th)=a\th'-1$, one can
easily check that $\a$ is compatible with the forms $\chi_{\th}$ and $\chi_{\th'}$ 
and hence with
maps $\phi_{\th}$ and $\phi_{\th'}$. This implies that $\a$ maps the division 
process of $[0,a-\th]$
associated with $\th$ to the division process of $[0,1]$ associated with $\th'$ as 
we claimed.  
By induction assumption $v'=m'\th'+n'$ where $m'<0$. Since $v'>0$ this implies 
that $n'>-m\th'>0$. Hence,
$$v=(a-\th)(m'\th'+n')=m'+(a-\th)n'=-n'\th+(m'+an')$$
has negative coefficient with $\th$.

Now let us prove that conversely every vector $v=m\th+n\in L_{\th}\cap (0,1)$ with 
$m<0$ belongs to $\BB_{\th}$.
We use induction in $|m|$. If $m=-1$ then $v$ coincides with the first division 
point $a-\th$, so the base
of induction is valid. Assume that the assertion is true for all vectors with 
smaller $|m|$ (and all $\th$).
Changing $\th$ to $-\th$ and $v$ to $1-v$
if necessary we can assume that $v<a-\th$ (note that $|m|$ remains invariant under
such a change). Now let us consider
$v'=v/(a-\th)\in L_{\th'}$ where $\th'=1/(a-\th)$. Then 
$v'=(ma+n)\th'-m$. We claim that 
$$m<ma+n<0.$$
This would finish the proof by applying the induction assumption to $v'$.
Since $m<0$ the inequalities we need are equivalent to
\begin{equation}
1>a+\frac{n}{m}>0.
\end{equation}
Now the inequalities $0<m\th+n<1$, $0<a-\th<1$ imply 
$$a+\frac{n}{m}<\th+1+\frac{n}{m}<1,$$
$$a+\frac{n}{m}>\th+\frac{n}{m}>\frac{1}{m}.$$
It remains to exclude the possibility $a+\frac{n}{m}=0$.
But in this case we would have $m\th+n=-m(a-\th)$ which contradicts
to $m\th+n<a-\th$.
\ed

\begin{cor} The set $\BB_{\th}$ is dense in $(0,1)$.
\end{cor}

\subsection{Construction of subbundles}
\label{subbun-sec}

We use two methods to construct subbundles in holomorphic vector bundles 
on $T=T_{\th,\tau}$. The first method is based on the binary division process described in the previous
section. 

\begin{thm}\label{main-thm} 
For every stable vector bundle $P$ on $T$ and every $r\in L_{\th}$ such that
$0<r<\rk P$ and $\chi(r,\rk P)>0$, there exists a subbundle $V\sub P$
such that $\rk V=r$.
\end{thm}

\Pf . Using Morita equivalences we can reduce ourselves to the case $\rk P=1$.
Then the condition $\chi(r,\rk P)>0$ on
$r=m\th+n\in\Z\th+\Z$ is equivalent to $m<0$. Thus, we have to show that for every such
$r<1$ there exists a subbundle $V\sub P$ with $\rk V=r$.
By Theorem \ref{main-div-thm} we have $r\in\BB_{\th}$.
Now we claim that one can associate to every division subsegment $[a,b]$ 
a stable vector bundle $V_{a,b}$ of rank $b-a$, such that $V_{0,1}=P$ and
for the new division point $c\in (a,b)$ one has an exact sequence
$$0\to V_{a,c}\to V_{a,b}\to V_{b,c}\to 0.$$
Indeed, assume that $V_{a,b}$ is already chosen. Then for the new division point $c\in (a,b)$
we can choose $V_{a,c}$ to be any stable bundle of rank $c-a$ and then observe that
the condition $\chi(c-a,b-a)=1$ implies that $\Hom(V_{a,c},V_{a,b})$ is one-dimensional
and $\Hom^1(V_{a,c},V_{a,b})=0$. We claim that the unique nonzero morphism
$f:V_{a,c}\to V_{a,b}$ is injective. Indeed, identifying $\Vect(T)$ with the subcategory
$\CC^{\th}\sub D^b(E)$ we can consider the cone of $f$ as an object in $D^b(E)$.
Since $\Cone(f)$ is the value on $V_{a,c}$ of the reflection functor associated with $V_{a,b}$
(see \cite{ST}), it follows that $\Cone(f)$ is stable. According to Lemma \ref{triangle-lem}
this implies that $\Cone(f)\in\CC^{\th}$ which implies our claim.
Now we can define inductively a family of subbundles
$V_{0,a}\sub V_{0,1}=P$ for all $a\in\BB_{\th}$, such that for every division subsegment $[a,b]$
one has $V_{0,a}\sub V_{0,b}$ and $V_{0,b}/V_{0,a}\simeq V_{a,b}$.
Indeed, assume that $V_{0,a}\sub V_{0,b}$ for the subsegment $[a,b]$ are already defined.
Then for the new division point
$c\in (a,b)$ we define $V_{0,c}\sub V_{0,b}$ as the preimage of $V_{a,c}\sub V_{a,b}$ under the projection $V_{0,b}\to V_{0,b}/V_{0,a}\simeq V_{a,b}$.
By the construction the subbundle $V=V_{0,r}$ has rank $r$.
\ed

Another way to construct subbundles
is based on the following lemma.

\begin{lem}\label{ss-lem} 
Let $A$ and $B$ be a pair of stable objects in $D^b(E)$.
Assume that
$\Hom(A,B)\neq 0$ and $\Hom^i(A,B)=0$ for $i\neq 0$. Then
for a generic morphism $f\in\Hom(A,B)$ the object $\Cone(f)$ is semistable.
\end{lem}

\Pf . Without loss of generality we can assume that $B=\OO_x$ for some
$x\in E$ (one has to use the action of a central extension of $\SL_2(\Z)$
on $D^b(E)$). Then $A$ is a stable vector bundle on $E$ and we have to
prove that for a generic morphism $f:A\to\OO_x$ the kernel of $f$ is
semistable. Let $r=\rk A$, $d=\deg A$. For every pair of relative prime
numbers $(r',d')$ let $M_{r',d'}$ be the moduli space
of stable bundles of rank $r'$ and
degree $d'$. Note that if $\ker(f)$ is unstable then it can be destabilized
by some stable bundle $V\in M_{r',d'}$ such that
$$\frac{d-1}{r}<\frac{d'}{r'}<\frac{d}{r}$$
and $r'<r$. Note that 
$$\dim\P\Hom(V,A)= dr'-rd'-1<r'-1.$$
Since $\dim M_{r',d'}=1$
we see that the family of possible bundles $V\sub A$ of this type has
dimension $<r'$. To every such $V$ there corresponds the $(r-r')$-dimensional
subspace in $\Hom(A,\OO_x)$ consisting of morphisms vanishing on $V$.
Therefore, a generic element $f\in\Hom(A,\OO_x)$ does not belong to any of
these subspaces. Hence, for such an element the bundle $\ker(f)$ will be
semistable. 
\ed

\begin{rem} The statement of the above lemma will become false if we
assume only that $A$ and $B$ are semistable. For example, any morphism
from $A=\OO_E$ to $B=\OO_x\oplus\OO_x$ has nonzero kernel and cokernel.
Hence, its cone is not semistable.
\end{rem}

The following result will play a crucial role in the proof of 
Theorem \ref{proj-thm}.

\begin{lem}\label{eps-lem} 
For every vector bundle $P$ and for a pair of real numbers 
$\eps>0$ and $C$ there exists a vector bundle $P'\sub P$
such that $\rk P'>\rk P-\eps$ and $P'$ is a direct sum of semistable
vector bundles of slopes $<C$.
\end{lem}

\Pf . First, 
we can reduce the proof to the case when the vector bundle $P$ is stable.
Indeed, it suffices to check that if $P$ fits into the exact sequence
$$0\to P_1\to P\to P_2\to 0$$
and the assertion of the lemma holds for $P_1$ and $P_2$ (and arbitrary $\eps$ and
$C$), then the assertion holds also for $P$. Let $A$ be the minimum of slopes of
semistable bundles in the Harder-Narasimhan filtration of $P_1$.
By assumption, we can find subbundles
$P'_1\sub P_1$ and $P'_2\sub P_2$ such that $\rk P'_1>\rk P_1-\eps/2$,
$\rk P'_2>\rk P_2-\eps/2$ and such that $P'_1$ and $P'_2$ are direct sums
of semistable bundles of slopes $<\min(A,C)$. Then 
$\Ext^1(P'_2,P_1)=0$, so there exists a splitting $P'_2\to P$.
Now we set $P'=P'_1\oplus P'_2\sub P$.

Thus, we can assume that $P$ is a stable holomorphic vector bundle.
Using principal convergents to $-\th$ we can choose a sequence of pairs of relatively prime integers
$(p_n,q_n)$ such that $p_n+q_n\th>0$, $\lim_{n\to\infty}(p_n+q_n\th)=0$
and $\lim_{n\to\infty}q_n=+\infty$. Let $V_n$ be a stable holomorphic vector bundle on $T$ with
$\rk(V_n)=p_n+q_n\th$. Then $\lim_{n\to\infty}\mu(V_n)=+\infty$. Hence,
for sufficiently large $n$ we have $\mu(V_n)>\mu(P)$ and $\rk(V_n)<\rk(P)$.
Now Lemmas \ref{ss-lem}
and \ref{triangle-lem} imply that a generic morphism $f:P\to V_n$ is surjective
and $P'=\ker(f)$ is semistable (provided $n$ is large enough). 
Note that $\deg(P')=\deg(P)-q_n\to-\infty$ as $n\to\infty$.
Hence, for large enough $n$ we will have 
$\mu(P')<C$ and $\rk(P')>\rk P-\eps$.
\ed

\subsection{Deformation equivalence}\label{deform-sec}

We start with the following combinatorial

\begin{lem}\label{comb-lem} 
Let us consider the set of unordered $n$-tuples $(v_1,\ldots,v_n)$
of primitive vectors in 
$\Z^2$ such that $0\not\in\R_{>0}v_1+\ldots+\R_{>0}v_n$
(we take the union over all $n\ge 1$). 
Consider the equivalence relation on this set generated by all 
relations of the following kind: an $n$-tuple $(v_1,\ldots,v_n)$ is equivalent
to $(w,\ldots,w,v_3,\ldots,v_n)$, where $w$ is repeated $m$ times, if
$v_1+v_2=mw$. Then $(v_1,\ldots,v_n)$ is equivalent to $(w_1,\ldots,w_m)$ iff
$v_1+\ldots+v_n=w_1+\ldots+w_m$.
\end{lem}

\Pf . We are going to show that every unordered $n$-tuple is equivalent
to an $n$-tuple of the form $(w,\ldots,w)$. Let us associate to every
$n$-tuple $(v_1,\ldots,v_n)$ a nonnegative integer by the following rule:
$$D(v_1,\ldots,v_n)=\sum_{i,j}\de_{\ge 0}(\det(v_i,v_j)),$$
where $\de_{\ge 0}(x)=x$ for $x\ge 0$ and $\de_{\ge 0}(x)=0$ for $x<0$.
It is clear that $D(v_1,\ldots,v_n)=0$ iff all $v_i$'s are the same.
Thus, it is enough to show that every $n$-tuple as above with not all $v_i$'s 
equal
is equivalent to some $m$-tuple $(w_1,\ldots,w_m)$ with
$D(w_1,\ldots,w_m)<D(v_1,\ldots,v_n)$.
Without loss of generality we can assume that $v_1\neq v_2$.
Since $v_1\neq -v_2$ we have $v_1+v_2=mw$ for some primitive vector $w$ and some
$m>0$. Consider the corresponding $n-2+m$-tuple $(w,\ldots,w,v_3,\ldots,v_n)$ 
equivalent
to $(v_1,\ldots,v_n)$. We claim that
$$D(w,\ldots,w,v_3,\ldots,v_n)<D(v_1,\ldots,v_n).$$
Indeed, since $D(v_1,v_2)>0$ it suffices to show that for every $i\ge 3$ one has
\begin{equation}\label{comb-ineq}
\begin{array}{l}
m[\de_{\ge 0}(\det(v_i,w))+\de_{\ge 0}(\det(w,v_i))]\le\nonumber\\
\de_{\ge 0}(\det(v_i,v_1))+\de_{\ge 0}(\det(v_i,v_2))
\de_{\ge 0}(\det(v_1,v_i))+\de_{\ge 0}(\det(v_2,v_i)).
\end{array}
\end{equation}
Assume for example that $\det(v_i,w)\ge 0$. Then
$\det(v_i,v_1+v_2)=m\det(v_i,w)\ge 0,$
so switching $v_1$ and $v_2$ if necessary we can assume that $\det(v_i,v_1)\ge 0$.
If $\det(v_i,v_2)\ge 0$ then \eqref{comb-ineq} holds since it fact it becomes
an equality. Finally, if $\det(v_i,v_2)\le 0$ then we have
$$m\det(v_i,w)=\det(v_i,v_1)+\det(v_i,v_2)\le\det(v_i,v_1)$$
which implies \eqref{comb-ineq}.
\ed

\begin{defi}
Let us say that two vector bundles on $T$ are {\it deformation equivalent}
if they belong to the same class with respect to the minimal equivalence relation
on isomorphism classes of vector bundles containing the following relations:

\noindent
(i) if $V_1$ and $V_2$ are stable and $\rk V_1=\rk V_2$ then $V_1\sim V_2$;

\noindent
(ii) if $0\to U\to V\to W\to 0$ is an exact triple of vector bundles then
$V\sim U\oplus W$;

\noindent
(iii) if $V_1\sim V_2$ then $U\oplus V_1\sim U\oplus V_2$ for any vector bundle
$U$.
\end{defi}

It is clear that deformation equivalent vector bundles have the same rank.
The following theorem states that the converse is also true.

\begin{thm}\label{def-thm} 
Let $V_1$ and $V_2$ be vector bundles on $T$ such that 
$\rk V_1=\rk V_2$. Then $V_1\sim V_2$.
\end{thm}

\Pf . The idea is to mimick the proof of Lemma \ref{comb-lem}.
It suffices to show that every vector bundle is deformation equivalent
to a bundle of the form $W_1\oplus\ldots \oplus W_m$
where $W_i$'s are stable and $\rk W_1=\ldots=\rk W_m$.
Using Harder-Narasimhan filtration and property (ii)
of our equivalence we can assume that our vector bundle is a direct sum
$V_1\oplus\ldots\oplus V_n$, where $V_i$'s are stable. Assume that
not all of them have the same rank, say, $\rk V_1\neq \rk V_2$.
Then we can reorder $V_1$ and $V_2$ in such a way that 
$\Hom(V_1,V_2)=0$ and $\Ext^1(V_1,V_2)\neq 0$.
According to Lemma \ref{ss-lem} for a generic extension
$$0\to V_2\to W\to V_1\to 0$$
the bundle $W$ is semistable. Using property (ii) we see that
$W\sim W_1\oplus\ldots\oplus W_m$, where $W_i$'s are stable and
have the same rank $(\rk V_1+\rk V_2)/m$. Thus, we have an equivalence
$$V_1\oplus\ldots\oplus V_n\sim W_1\oplus\ldots\oplus W_m\oplus V_3\oplus\ldots
\oplus V_n.$$
As we have seen in the proof of Lemma \ref{comb-lem}
repeating this procedure we will eventually arrive at the direct sum of
stable bundles of the same rank.
\ed

\section{Quasicoherent sheaves on $T$}
\label{quasi-sec}

\subsection{Ind-objects}\label{ind-sec}

Let us present some facts about ind-objects following sec.~8 of \cite{SGA4}.
Below by an inductive limit we always mean a small filtering inductive limit.
By a union of subobjects of an object in an abelian category we always mean a small
filtering union. 

Recall that with every category $\AA$ one can associate
the category $\Ind(\AA)$ of ind-objects of $\AA$ such that in $\Ind(\AA)$
all inductive limits exist (see \cite{SGA4}, sec.~8.2 and Prop.~8.5.1). 
The category
$\AA$ can be identified with a full subcategory of $\Ind(\AA)$.
Furthermore, if $\AA$ is abelian then so is $\Ind(\AA)$ (\cite{SGA4}, 8.9.9(c)) and the
natural embedding functor $\AA\to\Ind(\AA)$ is exact (\cite{SGA4}, Prop. 8.9.5(a)).
Also, Proposition~8.5.1 of \cite{SGA4} implies that  
if $(C_i)$ is an inductive system in $\Ind(\AA)$ then for $A\in\AA$ one has
\begin{equation}\label{lim-eq}
\dlim\Hom(A,C_i)\simeq\Hom(A,\dlim C_i)
\end{equation}

Assume that we are given a functor $f:\AA\to\CC$, where
$\CC$ is a category in which all inductive limits
exist. Then by \cite{SGA4}, 8.7.2, 
this functor extends to a functor $\ov{f}:\Ind(\AA)\to\CC$ commuting with inductive limits.
We will need the following result. 

\begin{prop}\label{ind-prop} 
(a) The functor $\ov{f}$ is fully faithful iff for every inductive system
$(A_i)$ in $\AA$  the natural map
$$\dlim\Hom_{\AA}(A,A_i)\to \Hom_{\CC}(f(A),\dlim f(A_i))$$
is bijective.

\noindent
(b) Assume that $\AA$ is abelian and the functor $f$ is exact. Then the functor
$\ov{f}$ is also exact.
\end{prop}

\Pf . (a) This is Prop. 8.7.5(a) of \cite{SGA4}.

(b) Since the embedding of $\AA$ into $\Ind(\AA)$ (resp., of $\CC$ into 
$\Ind(\CC)$) is exact, it is enough to check the exactness of
$\Ind(f):\Ind(\AA)\to\Ind(\CC)$. It remains to apply Cor. 8.9.8 of \cite{SGA4}.
\ed

Now let us consider the situation when we have an abelian category $\CC$ in which all
small direct sums (and hence, all small filtering inductive limits) exist. Let $\AA\sub\CC$
be a full abelian subcategory such that equivalent conditions of Proposition \ref{ind-prop}(a)
are satisfied. Then we can identify $\Ind(\AA)$ with a full subcategory of $\CC$.
In this situation we can give a convenient characterization of objects of $\CC$ that belong to
$\Ind(\AA)$ (see Propositions \ref{adm-prop} and \ref{fingen-adm-prop} below).

\begin{lem}\label{fingen-adm-lem} 
Let $A$ be an object of $\AA$, and let $S\sub A$ be a subobject in $\CC$. Then 
$A/S\in\wt{\AA}$ (equivalently, $S\in\wt{\AA}$)
iff $S=\cup_{i\in I}S_i$, where $(S_i)_{i\in I}$
is a filtering set of subobjects of $A$ in $\AA$.
\end{lem}

\Pf . First of all, we observe that since $\wt{\AA}\sub\CC$ is stable under kernels and cokernels,
we have $A/S\in\wt{\AA}$ iff $S\in\wt{\AA}$. If $S$ is a union
of subobjects $(S_i)$ of $A$, where $S_i\in\AA$ then $S\in\wt{\AA}$.
(because such a union can be taken in the category $\wt{\AA}$ and the embedding
$\wt{\AA}\to\CC$ is exact and commutes with inductive limits).
Conversely, assume that $S\in\wt{\AA}$. Then
$S=\dlim A_i$ for some inductive system $(A_i)$ in $\AA$.
Let $S_i=\im(A_i\to S)=\im(A_i\to A)$. Then $S_i\in\AA$ since $\AA$ is an abelian subcategory.
Also, $S=\cup_{i\in I}S_i$ (this union can be taken either in $\wt{\AA}$ or in $\CC$).
\ed

\begin{defi} Let us say that an object $C\in\CC$ is $\AA$-generated if there exists a surjection
$A\to C$ in $\CC$ with $A\in\AA$.
\end{defi}

\begin{prop}\label{fingen-adm-prop} 
In the above situation for an $\AA$-generated object $C\in\CC$ the following conditions are equivalent:

\noindent
(a) $C\in\wt{\AA}$;

\noindent
(b) $C\simeq A/S$, where $A\in\AA$ and $S$ is a union
of subobjects $(A_i)$ of $A$ with $A_i\in\AA$;

\noindent
(c) For every morphism $f: A\to C$ where $A\in\AA$, the kernel
$\ker(f)\sub A$ is a union of subobjects $(A_i)$ of $A$ with $A_i\in\AA$.
\end{prop}

\Pf . The implication (c)$\implies$(b) is clear, while
(b)$\implies$(a) follows from Lemma \ref{fingen-adm-lem}.
To show that (a)$\implies$(c) we note that if $C\in\wt{\AA}$ then $\ker(f)\in\wt{\AA}$. It remains
to apply Lemma \ref{fingen-adm-lem} again. 
\ed

\begin{prop}\label{adm-prop} 
For an object $C\in\CC$ the following conditions are equivalent:

\noindent
(a) $C\in\wt{\AA}$;

\noindent
(b) $C$ is a union $C=\cup_i C_i$, where
$C_i\in\wt{\AA}$ and $C_i$ is $\AA$-generated for all $i$;

\noindent
(c) $C$ is a union of a filtering set of $\AA$-generated subobjects,
and every $\AA$-generated subobject of $C$ belongs to $\wt{\AA}$;

\noindent
($c'$) $C$ is a union of a filtering set of $\AA$-generated subobjects and for every
morphism $f:A\to C$, where $A\in\AA$, the kernel $\ker(f)\sub A$ is a union of
subobjects $(A_i)$ of $A$ with $A_i\in\AA$.
\end{prop}

\Pf . The equivalence of (c) and (${\rm c'}$) follows from Proposition \ref{fingen-adm-prop}.
The implications (c)$\implies$(b) and (b)$\implies$(a) are clear. To prove (a)$\implies$(c)
we note that if $A\to C$ is a morphism, where $A\in\AA$ and $C\in\wt{\AA}$,
then its image also belongs to $\wt{\AA}$ (since the subcategory $\wt{\AA}\sub\CC$ is abelian). Hence, 
every $\AA$-generated subobject of $C$ belongs to $\wt{\AA}$. Also, if 
$C=\dlim A_i$, where $(A_i)_{i\in I}$ is an inductive system
in $\AA$, then $C=\cup_{i\in I} C_i$, where 
$C_i$ is the image of the morphism $A_i\to C$, so that each $C_i$ is $\AA$-generated.
\ed

\begin{ex} It is well known (and follows easily from Proposition \ref{adm-prop})
that the category of ind-coherent sheaves (ind-objects in $\Coh(X)$) on a Noetherian scheme $X$ is equivalent to the category $\Qcoh(X)$ of quasicoherent sheaves on $X$. Note that in this case the only $\Coh(X)$-generated objects in $\Qcoh(X)$ are coherent sheaves.
\end{ex}

\subsection{Holomorphic modules and holomorphic bundles on $T$}

To realize concretely ind-objects of the category of holomorphic vector bundles
on a noncommutative complex torus $T=T_{\th,\tau}$
we introduce an auxiliary category $HM(T)$ of {\it holomorphic
modules}. 

\begin{defi} 
An object of $HM(T)$ is a right 
$A_{\th}$-module $M$ equipped with a holomorphic structure, i.e., with
a map $\ov{\nabla}:M\to M$ satisfying the Leibnitz rule \eqref{Leibnitz-eq}.
The morphisms in $HM(T)$ are morphisms of $A_{\th}$-modules compatible
with holomorphic structures. 
\end{defi}

It is easy to see that $HM(T)$ is an abelian
category. By definition, a {\it holomorphic vector
bundle} on $T$ is a holomorphic
module $(M,\ov{\nabla})$ such that $M$ is a finitely generated projective
(right) $A_{\th}$-module. Thus, $\Vect(T)$ is a full subcategory in $HM(T)$. 
Also, it is easy to see that the natural embedding functor $\Vect(T)\to HM(T)$ is exact.

\begin{lem}\label{holom-map-lem}
Let $(M,\ov{\nabla})$ be a holomorphic module.
For every surjection of $A_{\th}$-modules
$\phi:A_{\th}^{\oplus n}\to M$ there exists a holomorphic structure
$\ov{\nabla}_{\phi}$ on $A_{\th}^{\oplus n}$ with respect to which
$\phi$ is holomorphic.
\end{lem}

\Pf .
Let $f_i=\phi(e_i)\in M$, where $e_1,\ldots,e_n$ is the standard basis of
$A_{\th}^{\oplus n}$. Then we have
$$\ov{\nabla}(f_j)=\sum_{i=1}^n f_i a_{ij}$$
for some $A_{\th}$-valued $n\times n$ matrix $M=(a_{ij})$.
Let us set for $(x_i)\in A_{\th}^{\oplus n}$
$$\ov{\nabla}_{\phi}(x_i)=(\de(x_i)+\sum_{j=1}^n a_{ij}x_j).$$
One can immediately check that the morphism 
$\phi$ becomes holomorphic with respect to $\ov{\nabla}_{\phi}$ and
$\ov{\nabla}$. 
\ed

Let us define the category $\Coh(T)$ of {\it coherent sheaves} on $T$
as the full subcategory in $HM(T)$ consisting of holomorphic modules 
$(M,\ov{\nabla})$ such that the $A_{\th}$-module $M$ is finitely
presented.

\begin{thm}\label{coh-thm} 
One has $\Coh(T)=\Vect(T)$, i.e., every coherent sheaf
on $T$ is a holomorphic vector bundle.
\end{thm}

\Pf . Let $(M,\ov{\nabla})$ be a coherent sheaf. Consider a finite
presentation 
$$M=\coker(f:A_{\th}^{\oplus m}\to A_{\th}^{\oplus n}).$$
Applying Lemma \ref{holom-map-lem} to the natural projection
$A_{\th}^{\oplus n}\to M$ we find a holomorphic structure 
$\ov{\nabla}_2$ on $A_{\th}^{\oplus n}$ with respect to which this
projection is holomorphic. Now the submodule $\im(f)\sub A_{\th}^{\oplus n}$
is preserved by $\ov{\nabla}_2$, so we can view $(\im(f),\ov{\nabla}_2)$ 
as a holomorphic module. Applying Lemma \ref{holom-map-lem} to
the surjection $f:A_{\th}^{\oplus m}\to \im(f)$ we find a holomorphic
structure $\ov{\nabla}_1$ on $A_{\th}^{\oplus m}$ such that $f$
is holomorphic with respect to $\ov{\nabla}_1$ and $\ov{\nabla}_2$.
Hence, $f$ becomes a morphism in the category $\Vect(T)$.
Since $\Vect(T)$ is abelian, it follows that $M=\coker(f)$ is an object of
$\Vect(T)$.
\ed

\begin{cor}\label{fingen-proj-cor} Let $P$ be a holomorphic bundle
and let $S\sub P$ be a finitely generated holomorphic submodule.
Then $S$ is a direct summand of $P$ as an $A_{\th}$-module. Hence,
it is a holomorphic subbundle of $P$.
\end{cor}

\Pf . The holomorphic module $P/S$ is finitely presented, hence, it
is a vector bundle by the above theorem.
\ed

We equip the trivial module $A_{\th}$ 
with the standard holomorphic structure $\de_{\tau}$
Thus, a {\it holomorphic ideal} in $A_{\th}$ is a right ideal
$I\sub A_{\th}$ such that $\de_{\tau}(I)\sub I$.

\begin{cor}\label{fingen-cor} 
Let $I\sub A_{\th}$ be a finitely generated
holomorphic ideal. Then there exists an idempotent $e\in A_{\th}$ such that 
$I=e A_{\th}$.
\end{cor}

\begin{rems}
1. Let us say that an idempotent $e\in A_{\th}$ is (right) {\it holomorphic} if
$$e\de_{\tau}(e)=\de_{\tau}(e),$$
or equivalently,
$\de_{\tau}(e)\in e A_{\th}$. The above theorem shows that a map
$e\mapsto e A_{\th}$ gives a surjection from the set of holomorphic idempotents
to that of finitely generated holomorphic ideals.
It is easy to see that the fiber of this map over $e A_{\th}$ coincides
with $e+e A_{\th}(1-e)$.

\noindent
2. Proposition \ref{quasi-proj-mod-prop} below implies that countable
filtering unions of finitely generated holomorphic ideals are still projective $A_{\th}$-modules.
However, we will see that such ideals 
are not necessarily direct summands in $A_{\th}$
(see Theorem \ref{quasi-main-thm}).
\end{rems}

\subsection{Quasicoherent sheaves as holomorphic modules}\label{quasi-adm-sec}

We define the category of
{\it quasicoherent sheaves} on $T$ by setting $\Qcoh(T)=\Ind(\Vect(T))$.
Thus, by definition, quasicoherent sheaves on $T$ are ind-objects
in the category $\Vect(T)$. They form an abelian category containing
$\Vect(T)$ as a full subcategory.
We are going to give a more concrete realization of $\Qcoh(T)$ using holomorphic modules.

It is easy to see that in $HM(T)$
all small filtering inductive limits exist. More precisely,
the natural embedding of $HM(T)$ into the category of $A_{\th}$-modules
is exact and commutes with small inductive limits.
Therefore, the natural fully faithful exact
embedding $\Vect(T)\hra HM(T)$ extends to
an exact functor $\Qcoh(T)\to HM(T)$ (see Proposition \ref{ind-prop}(b)).

\begin{prop}\label{fully-faith-prop} 
The above functor $\Qcoh(T)\to HM(T)$ is fully faithful.
\end{prop}

\Pf . According to Proposition \ref{ind-prop}(a) we have to check that for every
small filtering inductive system $(P_i)_{i\in I}$ in $\Vect(T)$ and every
$P\in\Vect(T)$ the canonical map
\begin{equation}\label{limmap}
\dlim\Hom_{\Vect(T)}(P,P_i)\to\Hom_{HM(T)}(P,\dlim P_i)
\end{equation}
is an isomorphism. It is easy to see that if we replace morphisms
in $HM(T)$ by morphisms of $A_{\th}$-modules then the similar map
is an isomorphism because $P$ is a finitely generated projective
$A_{\th}$-module. This immediately implies injectivity of \eqref{limmap}.
To check surjectivity let us assume that 
$f:P\to\dlim P_i$ is any morphism in $HM(T)$. Then we can lift $f$ to
a morphism of $A_{\th}$-modules $f_{i_0}:P\to P_{i_0}$. Let $e_1,\ldots,e_n$
be generators of $P$. Then we have
$$\ov{\nabla}(e_i)=\sum_j e_j a_{ij}$$
for some $a_{ij}\in A_{\th}$. Hence,
$$\ov{\nabla}(f(e_i))=\sum_j f(e_j) a_{ij}.$$
It follows that for some $i_1>i_0$ we will have
$$\ov{\nabla}(f_{i_1}(e_i))=\sum_j f_{i_1}(e_j) a_{ij},$$
where $f_{i_1}:P\to P_{i_1}$ is the morphism of $A_{\th}$-modules induced
by $f_{i_0}$. But this means that $f_{i_1}$ is compatible with holomorphic
structures, so it is a morphism in $HM(T)$.
\ed

\begin{defi}
Let us call a holomorphic module {\it admissible} if
it can be represented as an inductive limit of holomorphic bundles.
\end{defi}

\begin{cor}\label{adm-cor}
The category $\Qcoh(T)$ is equivalent
to the full subcategory of admissible modules in $HM(T)$.
\end{cor}

Let us say that a holomorphic module $M$ is {\it finitely generated} if it is
finitely generated as an $A_{\th}$-module. The following lemma shows that
this is equivalent to $M$ being $\Vect(T)$-generated.

\begin{lem}\label{fingen-mod-lem} 
A holomorphic module $M$ is finitely generated iff there exists
a holomorphic bundle $P$ and a surjection $P\to M$ in $HM(T)$.
\end{lem}

\Pf . The ``if'' part is clear. The ``only if'' part follows immediately
from Lemma \ref{holom-map-lem}.
\ed

Applying Propositions \ref{fingen-adm-prop} and \ref{adm-prop} we get 
the following characterization of admissible
holomorphic modules. 

\begin{prop}\label{fingen-adm-prop2} 
For a finitely generated holomorphic module $M$ the following conditions are equivalent:

\noindent
(a) $M$ is admissible;

\noindent
(b) $M\simeq P/S$, where $P$ is a holomorphic vector bundle, and $S$ is a union
of holomorphic subbundles in $P$;

\noindent
(c) For every morphism $f: P\to M$, where $P$ is a holomorphic vector bundle, 
$\ker(f)\sub P$ is a union of holomorphic subbundles in $P$.
\end{prop}

\begin{prop}\label{adm-prop2} 
For a holomorphic module $M$ the following conditions are equivalent:

\noindent
(a) $M$ is admissible;

\noindent
(b) $M$ is the union of its finitely generated admissible holomorphic submodules.

\noindent
(c) $M$ is the union of its finitely generated holomorphic submodules
and every such submodule is admissible;

\noindent
($c'$) $M$ is the union of its finitely generated holomorphic submodules and for every
morphism $f:P\to M$, where $P$ is a holomorphic vector bundle, $\ker(f)$ is a union of
holomorphic subbundles in $P$.
\end{prop}

\begin{cor} Let $0\to M'\to M\to M"\to 0$ be an exact sequence of holomorphic modules
such that $M'$ and $M"$ are admissible. Then $M$ is also admissible.
\end{cor}

\Pf . First, let us show that every element $x\in M$ is contained in
a finitely generated holomorphic submodule. Let $N"\sub M"$ be
a finitely generated holomorphic submodule containing the image of $x$ in $M"$.
We can lift generators of $N"$ to some elements $e_1,\ldots,e_n\in M$.
Then $x=\sum a_i e_i +y$ and
$\ov{\nabla}e_j=\sum a_{ij} e_i+y_j$ for some $a_i, a_{ij}\in A$ and
$y, y_j\in M'$.
Let $N'$ be a finitely generated holomorphic submodule in $M'$ containing $y$ and
$(y_j)$. Then the $A_{\th}$-submodule generated by $N'$ and $(e_i)$ is holomorphic and
contains $x$. Next, let $f: P\to M$ be any morphism, where $P$ is a holomorphic vector bundle.
We have to show that $\ker(f)$ is admissible. Let $f":P\to M"$ and $f':\ker(f")\to M'$ be
the induced morphism. Then $\ker(f")$ is admissible and hence 
$\ker(f')$ is admissible. It remains to observe that $\ker(f)=\ker(f')$.
\ed

It is not clear whether there exists a simple
characterization of all $A_{\th}$-modules
underlying quasicoherent sheaves. One obvious condition is flatness (since
inductive limits of projective modules are flat).
In section \ref{count-sec} we will introduce an abelian subcategory of
countably presented quasicoherent sheaves and will show
that projective dimension of $A_{\th}$-modules underlying such sheaves is $\le 1$
(see Theorem \ref{count-thm}).

\subsection{Derived categories equivalence}\label{FM-sec}

Recall that the category $\Vect(T)$ is equivalent to the abelian subcategory $\CC^{\th}$ in 
the derived category $D^b(E)$ of coherent sheaves on the elliptic curve $E=\C/(\Z+\Z\tau)$.
Moreover, we have an equivalence of derived categories 
\begin{equation}\label{derived-eq}
D^b(\Vect(T))\simeq D^b(E)
\end{equation}
(see Proposition 1.4 of \cite{P-real}).
In this section we establish a similar equivalence for the category of quasicoherent
sheaves $\Qcoh(T)$.

Let us consider a more general situation.
Let $\AA$ be a Noetherian abelian category and let
$\wt{\AA}=\Ind(\AA)$. We view $\AA$ as a full subcategory in $\wt{\AA}$ in a standard way.

\begin{lem}\label{Noeth-lem} 
(i) For every surjection $A\to B$ in $\wt{\AA}$ with $B\in\AA$ there exists a subobject
$A'\sub A$ such that $A'\in\AA$ and the induced map $A'\to B$ is surjective.

\noindent
(ii) The subcategory $\AA\sub\wt{\AA}$ is closed under passing to quotients and subobjects, and
under extensions.
\end{lem}

\Pf . (i) Let $A=\dlim A_i$, where $(A_i)$ is an inductive system in $\AA$. Let
$B_i\sub B$ be the image of the induced map $A_i\to A$. Then $(B_i)$ is a filtering
set of subobjects of $B$ such that $\cup_i B_i=B$. Hence, $B_i=B$ for some $i$,
so we can take as $A'$ the image of the corresponding morphism $A_i\to A$. 

\noindent
(ii) Lemma \ref{fingen-adm-lem} easily implies that $\AA$ is stable under quotients and subobjects in 
$\wt{\AA}$. Assume that we have an exact sequence 
$$0\to A\to \wt{B}\to C\to 0$$
in $\wt{\AA}$ with $A,C\in\AA$. By part (i) we can find a subobject $B\sub\wt{B}$ such that
the induced map $B\to C$ is surjective. Then $\wt{B}$ is isomorphic to a quotient of $A\oplus B$,
so $\wt{B}\in\AA$.
\ed

\begin{lem}\label{ext-lim-lem} 
Let $(X_i)$ be an inductive system in $\wt{\AA}$. Then for every $A\in\AA$ the natural
map
$$\dlim\Ext^1(A,X_i)\to\Ext^1(A,\dlim X_i)$$
is an isomorphism.
\end{lem}

\Pf . First, let us check surjectivity. Assume that we have an extension 
$$0\to\dlim X_i\to\wt{E}\to A\to 0$$
in $\wt{\AA}$. By Lemma \ref{Noeth-lem}(i) there exists a subobject $E\sub\wt{E}$ with
$E\in\AA$ that still surjects onto $A$. Let $S=E\cap\dlim X_i\sub E$. Then the above extension
is the push-out of an extension
\begin{equation}\label{extension-eq}
0\to S\to E\to A\to 0
\end{equation}
by the inclusion map $S\to\dlim X_i$. But the latter map factors through a map $S\to X_i$.
Taking the push-out of \eqref{extension-eq} by this map we get an element
of $\Ext^1(A,X_i)$ inducing the original extension.

To check injectivity let us consider an extension
\begin{equation}\label{extension-eq2}
0\to X_{i_0}\to \wt{E}\to A\to 0
\end{equation}
such that its push-out by $X_{i_0}\to\dlim X_i$ becomes trivial.
As above we can find an extension \eqref{extension-eq} 
such that \eqref{extension-eq2} is its push-out by a map $S\to X_{i_0}$.
The fact that the push-out by $f:S\to X_{i_0}\to\dlim X_i$ becomes trivial means
that there is a map $s: E\to\dlim X_i$ extending $f:S\to\dlim X_i$.
But $s$ factors through some map $E\to X_j$. Furthermore, choosing $X_j$
sufficiently far in the inductive system will guarantee that this map restricts to the map
$g:S\to X_{i_0}\to X_j$. Hence, the push-out of \eqref{extension-eq} by $g$ is trivial.
Thus, the induced element of $\Ext^1(A,X_j)$ is trivial.
\ed

Let $(\TT,\FF)$ be a torsion pair in $\AA$. Recall (see \cite{HRS}) that this means that
$\TT$ and $\FF$ are full subcategories in $\AA$ such that $\Hom(T,F)=0$ for $T\in\TT$, $F\in\FF$,
and every object $X\in\AA$ fits into an exact triangle 
$$0\to T\to X\to F\to 0$$
with $T\in\TT$, $F\in\FF$. 
With a torsion pair $(\TT,\FF)$ one associates the
{\it tilted} abelian subcategory $\AA_t\sub D^b(\AA)$ by setting 
$$\AA_t=\{ K\in D^b(\AA)\ |\ H^{-1}K\in\FF, H^0K\in\TT, H^i K=0\text{ for }i\neq -1,0\}.
$$ 
The pair of subcategories $(\FF[1],\TT)$ is a torsion pair in $\AA_t$.
A torsion pair $(\TT,\FF)$ is called {\it cotilting} (resp., {\it tilting}) if every object of $\AA$ is a quotient of an object in $\FF$ (resp., subobject of an object in $\TT$).
A torsion pair $(\TT,\FF)$ in $\AA$ is cotilting iff the pair $(\FF[1],\TT)$ in $\AA_t$ is tilting
(see \cite{HRS}, Prop. I.3.2).

Set 
$$\wt{\TT}=\{ X\in\wt{\AA}\ |\ \text{ there exists a surjection } \oplus_i T_i\to X, \text{ with } T_i\in\TT\text{ for all }i\},$$
and let $\wt{\FF}\sub\wt{\AA}$ be the right orthogonal of $\TT$ in $\wt{\AA}$.
 
\begin{lem}\label{torsion-lem} 
(i) The subcategories $\wt{\TT}$ and $\wt{\FF}$ are closed under inductive limits and
the natural functors $\Ind(\TT)\to\wt{\TT}$, $\Ind(\FF)\to\wt{\FF}$ are equivalences.

\noindent
(ii) $(\wt{\TT},\wt{\FF})$ is a torsion pair in $\wt{\AA}$. We have $\wt{\TT}\cap\AA=\TT$,
$\wt{\FF}\cap\AA=\FF$. Furthermore, if $(\TT,\FF)$ is cotilting then so is
$(\wt{\TT},\wt{\FF})$. 

\noindent
(iii) Let $\wt{\AA}_t\sub D^b(\wt{\AA})$ be the tilted abelian category associated with
$(\wt{\TT},\wt{\FF})$. Then in $\wt{\AA}_t$ arbitrary (small) direct sums exist and
for any set $I$ the direct sum functor $(A_i)_{i\in I}\to \oplus_{i\in I} A_i$ is exact.
In other words, $\wt{\AA}_t$ satisfies the axiom AB4 (see \cite{Gr}, 1.5).
\end{lem}

\Pf . (i) It is clear that $\wt{\TT}$ is closed under direct sums. Since we have a surjection
$\oplus_i X_i\to\dlim X_i$ for an inductive system $(X_i)$, it follows that $\wt{\TT}$
is closed under inductive limits. The same assertion for $\wt{\FF}$ follows from \eqref{lim-eq}.
The natural functors $\Ind(\TT)\to\wt{\AA}$ and $\Ind(\FF)\to\wt{\AA}$ are fully faithful
by Proposition \ref{ind-prop}(a). It remains to show that every object of $\wt{\TT}$
(resp., $\wt{\FF}$) can be represented as the limit of an inductive system $(T_i)$ with $T_i\in\TT$
(resp., $(F_i)$ with $F_i\in\FF$). If $X\in\wt{\TT}$ then we have a surjection
$\oplus_i T_i\to X$. Let $T'_i$ be the image of the map $T_i\to X$. Then $X=\dlim T'_i$ and
$T'_i\in\TT$ since $\TT$ is closed under quotients. On the other hand, if $X\in\wt{\FF}$ then
$X=\dlim_i X_i$, where $(X_i)$ is a filtering set of $\AA$-generated subobjects in $X$
(see Proposition \ref{adm-prop}). By Lemma \ref{Noeth-lem}(ii) we have $X_i\in\AA$.
Since $X_i$ is a subobject of $X\in\wt{\FF}$, it belongs to the right orthogonal of $\TT$ in $\AA$.
Hence, $X_i\in\FF$.

\noindent (ii) If $X\in\wt{\TT}$ and $Y\in\wt{\FF}$ then there exists a surjection
$\oplus_i T_i\to X$, where $T_i\in\TT$. Since $\Hom(T_i,Y)=0$ this implies that $\Hom(X,Y)=0$.
If $X\in\wt{\AA}$ is an arbitrary object then $X=\dlim X_i$ for an inductive system $(X_i)$ in $\AA$.
For every $i$ we have an exact sequence
$$0\to A_i\to X_i\to B_i\to 0$$
with $A_i\in\TT$, $B_i\in\FF$. Moreover, these objects assemble into inductive systems $(A_i)$
and $(B_i)$, and the sequence
$$0\to\dlim A_i\to\dlim X_i\to\dlim B_i\to 0$$
in $\wt{\AA}$ is still exact (the exactness on the left follows from \eqref{lim-eq}).
Since $\dlim A_i\in\wt{\TT}$ and $\dlim B_i\in\wt{\FF}$, it follows that $(\wt{\TT},\wt{\FF})$ is
a torsion pair in $\wt{\AA}$. 

It is clear that $\wt{\FF}\cap\AA=\FF$ and that $\TT\sub\wt{\TT}\cap\AA$.
On the other hand, $\wt{\TT}\cap\AA$ is left orthogonal to $\FF$, hence it is contained in $\TT$.

Finally, assume $(\TT,\FF)$ is cotilting. Then for every $X\in\wt{\AA}$ there exists an inductive system
$(A_i)$ in $\AA$ such that $X=\dlim A_i$. Pick a surjection $F_i\to A_i$ with $F_i\in\FF$
for every $i$. Then the composed map
$$\oplus F_i\to\oplus A_i\to X$$
is a surjection and $\oplus F_i\in\FF$.

\noindent (iii) 
Note that the category $\wt{\AA}=\Ind(\AA)$ 
satisfies the axiom AB4. Indeed, the functor of direct sum is always right exact. The fact that in this case it is left exact follows easily from \eqref{lim-eq}. 
Therefore, in the (unbounded) derived category $D(\wt{\AA})$
arbitrary (small) direct sums exist and direct sums
of triangles are triangles (see \cite{BN}, Cor.~1.7). Also, by construction these direct sums commute
with the cohomology functors. Using (i) we deduce that $\wt{\AA}_t\sub D(\wt{\AA})$ is stable
under direct sums.
Now let $X_i\hra Y_i$, $i\in I$, be a collection of injective morphisms in $\wt{\AA}_t$.
Then we have exact triangles of the form
$$X_i\to Y_i\to Z_i\to X_i[1]$$
with $Z_i\in\wt{\AA}_t$. Taking the direct sums over $i$ we get an exact triangle
$$\oplus_i X_i\to \oplus_i Y_i\to\oplus_i Z_i\to\oplus_i X_i[1].$$
Since $\oplus_i Z_i$ belongs to $\wt{\AA}_t$, this implies that $\oplus_i X_i\to\oplus_i Y_i$ is
injective.
\ed

\begin{thm}\label{FM-thm}
Let $\AA$ be a Noetherian abelian category
of homological dimension $\le 1$,
$(\TT,\FF)$ a cotilting torsion pair in $\AA$.
Set $\wt{\AA}=\Ind(\AA)$ and define the torsion pair $(\wt{\TT},\wt{\FF})$ in $\wt{\AA}$ as above.
Let $\AA_t$ (resp., $\wt{\AA}_t$) be the tilted abelian category associated with
$(\TT,\FF)$ (resp., $(\wt{\TT},\wt{\FF})$).
Then there is an equivalence of categories 
$$\wt{\AA}_t\simeq\Ind(\AA_t)$$
and an exact equivalence of triangulated categories
$$D^b(\wt{\AA})\simeq D^b(\Ind(\AA_t)).$$
\end{thm}

\Pf . We have the natural exact functor
$\AA_t\to\wt{\AA}_t$. Since $\wt{\AA}_t$ is closed under direct sums, it extends to
a functor $\Phi:\Ind(\AA_t)\to\wt{\AA}_t$. To check that this functor is fully faithful we have
to check that for $X\in\AA_t$ and an inductive system $(Y_i)$ in $\AA_t$ the natural map
\begin{equation}\label{tilt-lim-eq}
\dlim\Hom_{\AA_t}(X,Y_i)\to\Hom_{\wt{\AA}_t}(X,\dlim Y_i)
\end{equation}
is an isomorphism (see Proposition \ref{ind-prop}).
Note that the category $\AA_t$ (resp., $\wt{\AA}_t$) is equipped with a torsion pair
$(\FF[1],\TT)$ (resp., $(\wt{\FF}[1],\wt{\TT})$). 
Thus, for every $X\in\AA_t$ we have an exact sequence 
\begin{equation}\label{AXB}
0\to A\to X\to B\to 0
\end{equation}
in $\AA_t$ with $A\in\FF[1]$ and $B\in\TT$.

\noindent
{\bf Step 1}. {\it The map \eqref{tilt-lim-eq} is an isomorphism if all $Y_i$ belong to $\FF[1]$.}
For the proof consider the morphism of the long exact sequences associated with \eqref{AXB}: 
\begin{equation}
\begin{diagram}
0 \to &\dlim\Hom_{\AA_t}(B,Y_i)&\to &\dlim\Hom_{\AA_t}(X,Y_i)&\to&\dlim\Hom_{\AA_t}(A,Y_i)&\to&\dlim\Ext^1_{\AA_t}(B,Y_i)\\
&\dTo{\a} && \dTo && \dTo{\b} && \dTo{\ga}\\
0 \to &\Hom_{\wt{\AA}_t}(B,\dlim Y_i)&\to &\Hom_{\wt{\AA}_t}(X,\dlim Y_i)&\to&
\Hom_{\wt{\AA}_t}(A,\dlim Y_i)&\to&\Ext^1_{\wt{\AA}_t}(B,\dlim Y_i) \nonumber
\end{diagram}
\end{equation}
By the five-lemma it is enough to check that $\a$ and $\b$ are isomorphisms and that $\ga$ is
injective. The fact that $\b$ is an isomorphism is the consequence of the equivalence
$\wt{\FF}\simeq\Ind(\FF)$ (see Lemma \ref{torsion-lem}(i)). 
The fact that $\a$ is an isomorphism follows from Lemma
\ref{ext-lim-lem} along with the fact that $\Ext^1_{\AA}(B,Y_i)\simeq\Ext^1_{\wt{\AA}}(B,Y_i)$
since $\AA$ is closed under extensions in $\wt{\AA}$ (by Lemma \ref{Noeth-lem}(ii)). Finally,
$$\Ext^1_{\AA_t}(B,Y_i)\simeq\Ext^1_{D^b(\AA)}(B,Y_i),$$
so the vanishing of $\Ext^2$ in $\AA$ implies the vanishing of all these groups.

\noindent
{\bf Step 2}. {\it If $(\wt{Y}_i)$ is an inductive system in $\wt{\FF}[1]$ and $X\in\AA_t$ then the natural map
$$\dlim\Hom_{\wt{\AA}_t}(X,\wt{Y}_i)\to\Hom_{\wt{\AA}_t}(X,\dlim \wt{Y}_i)$$
is injective.} The proof is similar to Step 1 but easier: it is enough to consider only the first three terms
in the long exact sequences associated with \eqref{AXB} and to use Lemma \ref{ext-lim-lem}.

\noindent
{\bf Step 3}. {\it The map \eqref{tilt-lim-eq} is injective.}
Indeed, our assumption that $(\TT,\FF)$ is a cotilting torsion pair in $\AA$ implies that
$(\FF[1],\TT)$ is a tilting torsion pair in $\AA_t$, i.e., every object of $\AA_t$ can be embedded
into an object of $\FF[1]$. Let us choose  for every $i$
an embedding $Y_i\sub Y'_i$ with $Y'_i\in\FF[1]$.
We want to replace these embeddings by a new collection of embeddings $Y_i\sub\wt{Y}_i$ in 
$\wt{\AA}_t$, where $\wt{Y}_i\in\wt{\FF}[1]$, such that $(\wt{Y}_i)$ form an inductive system
and the above embeddings give a morphism of inductive systems. This is achieved by defining
$\wt{Y}_i$ from the following pushout diagram
\begin{equation}
\begin{diagram}
\oplus_{j\le i} Y_j & \rTo{} & Y_i\\
\dTo{} & &\dTo{}\\
\oplus_{j\le i} A_j & \rTo{} & \wt{Y}_i\nonumber
\end{diagram}
\end{equation}
Since the horizontal arrows in this diagram are surjective, we have $\wt{Y}_i\in\wt{\FF}[1]$
(because $\wt{\FF}[1]$ is the left orthogonal to $\wt{\TT}$ in $\wt{\AA}_t$).
On the other hand, by Lemma \ref{torsion-lem}(iii) the left vertical arrow is injective, hence
the right vertical arrow $Y_i\to\wt{Y}_i$ is also injective. It is easy to see that these maps
form a morphism of inductive systems.
We can include the map \eqref{tilt-lim-eq} into a commutative square
\begin{equation}
\begin{diagram}
\dlim\Hom_{\AA_t}(X,Y_i) & \rTo{} & \Hom_{\wt{\AA}_t}(X,\dlim Y_i)\\
\dTo{} & & \dTo{} \\
\dlim\Hom_{\wt{\AA}_t}(X,\wt{Y}_i) & \rTo{} & \Hom_{\wt{\AA}_t}(X,\dlim \wt{Y}_i)\nonumber
\end{diagram}
\end{equation}
Now the required injectivity follows immediately from the injectivity of the left vertical
arrow and from Step 2.

\noindent
{\bf Step 4}. {\it The map \eqref{tilt-lim-eq} is an isomorphism}.
For the proof let us consider for every $i$ the exact sequence
$$0\to A_i \to Y_i\to B_i\to 0$$
with $A_i\in\FF[1]$ and $B_i\in\TT$. The objects $(A_i)$ (resp., $(B_i)$)
form an inductive system and we claim that the limit sequence
\begin{equation}\label{lim-ex-seq}
0\to \dlim A_i\stackrel{f}{\to}\dlim Y_i\to\dlim B_i\to 0
\end{equation}
in $\wt{\AA}_t$ is exact. Indeed, the only problem is the injectivity of $f$.
By Step 1 we obtain that for every $X\in\AA_t$ one has
$$\dlim \Hom_{\AA_t}(X,A_i)\simeq\Hom_{\wt{\AA}_t}(X,\dlim A_i).$$
Together with Step 3 this implies that the morphism
$$\Hom_{\wt{\AA}_t}(X,\dlim A_i)\to\Hom_{\wt{\AA}_t}(X,\dlim Y_i)$$
induced by $f$ is injective. Therefore, $\Hom_{\wt{\AA}_t}(X,\ker(f))=0$
for every $X\in\AA_t$. Consider the exact sequence
$$0\to C\to\ker(f)\to D\to 0$$
in $\wt{\AA}_t$ with $C\in\wt{\FF}[1]$ and $D\in\wt{\TT}$. Since
$\wt{\FF}\simeq\Ind(\FF)$ by Lemma \ref{torsion-lem}(i), it follows that $C=0$, i.e., $\ker(f)\in\wt{\TT}$.
Similarly, using the fact that $\wt{\TT}\simeq\Ind(\TT)$ we conclude that $\ker(f)=0$.
Next, let us consider the morphism of the long exact sequences
\begin{equation}
\begin{diagram}
0 \to &\dlim\Hom_{\AA_t}(X,A_i)&\to &\dlim\Hom_{\AA_t}(X,Y_i)&\to&\dlim\Hom_{\AA_t}(X,B_i)&\to&\dlim\Ext^1_{\AA_t}(X,A_i)\\
&\dTo{} && \dTo && \dTo{\de} && \dTo{}\\ 
0 \to &\Hom_{\wt{\AA}_t}(X,\dlim A_i)&\to &\Hom_{\wt{\AA}_t}(X,\dlim Y_i)&\to&
\Hom_{\wt{\AA}_t}(X,\dlim B_i)&\to&\Ext^1_{\wt{\AA}_t}(X,\dlim A_i)\nonumber
\end{diagram}
\end{equation}
By the five-lemma and Steps 1 and 3 it remains to check that $\de$ is an isomorphism.
Consider the exact sequence \eqref{AXB} again.
Then we have
$$\dlim\Hom_{\AA_t}(X,B_i)\simeq\dlim\Hom_{\TT}(B,B_i)\simeq\Hom_{\wt{\TT}}(B,\dlim B_i)\simeq
\Hom_{\wt{\AA}_t}(X,\dlim B_i),$$
where we used again the fact that $\wt{\TT}\simeq\Ind(\TT)$.

\noindent
{\bf Step 5}. {\it The functor $\Phi:\Ind(\AA_t)\to\wt{\AA}_t$ is essentially surjective on objects}.
Indeed, since the torsion pair $(\TT,\FF)$ is cotilting the same is true about the torsion pair
$(\wt{\TT},\wt{\FF})$ in $\wt{\AA}$ (see Lemma \ref{torsion-lem}(ii)). 
Therefore, the torsion pair $(\wt{\FF}[1],\wt{\TT})$ in 
$\wt{\AA}_t$ is tilting. So for any object $X\in\wt{\AA}_t$ we can find an embedding
$X\sub A$ with $A\in\wt{\FF}[1]$. Since $\wt{\FF}[1]$ is closed under quotients we obtain
that $X$ is the kernel of a morphism $f:A\to A'$ with $A,A'\in\wt{\FF}[1]$. Note that
$\wt{\FF}[1]$ is contained in the essential image of $\Phi$
(since $\wt{\FF}\simeq\Ind(\FF)$ by Lemma \ref{torsion-lem}(i)). Since by Step 4 the functor
$\Phi$ is fully faithful, we derive that the morphism $f$ is contained in the image of $\Phi$. Finally,
exactness of $\Phi$ (see Proposition \ref{ind-prop}(b))
implies that $X\simeq\ker(f)$ is also contained in the image of $\Phi$.

The statement about the equivalence of derived categories is an immediate consequence
of Proposition 5.4.3 of \cite{BVdB}.
\ed

The category $\Vect(T)$ is equivalent to the tilted abelian category associated with the cotilting torsion
pair $(\Coh_{>\th}(E),\Coh_{<\th}(E))$ in the category $\Coh(E)$ of the coherent sheaves on the
elliptic curve $E$,
where $\Coh_{<\th}(E)$ consists of direct sums of semistable bundles of slopes $<\th$ and
$\Coh_{>\th}(E)$ consists of direct sums of torsion sheaves and semistable bundles of slopes $>\th$
(see \cite{P-holbun} or \cite{P-real}, sec.~1.2).
Let $\Qcoh(E)\simeq\Ind(\Coh(E))$ be the category of quasicoherent sheaves on $E$.
It is easy to see that the assumptions of Theorem \ref{FM-thm} are satisfied, so we derive
the following result.

\begin{thm}\label{FM-tori-thm} 
The category $\Qcoh(T)$ is equivalent to the tilted abelian category associated with
the torsion pair $(\Qcoh_{>\th}(E),\Qcoh_{<\th}(E))$ in $\Qcoh(E)$, where
$\Qcoh_{<\th}(E)$ consists of torsion free quasicoherent sheaves $F$, such that any semistable subsheaf of $F$ has slope $<\th$, and 
$\Qcoh_{>\th}(E)$ consists of quotients of arbitrary (small) direct sums of torsion sheaves and
semistable bundles of slopes $>\th$. 
There is also an exact equivalence of derived categories 
$$D^b(\Qcoh(T))\simeq D^b(\Qcoh(E))$$
extending \eqref{derived-eq}.
\end{thm}

\subsection{The rank function}\label{rank-sec}

The main result of this section is the following

\begin{thm}\label{rank-thm} 
There exists a unique extension of the function
$\rk$ on holomorphic bundles over $T$ to 
a function on quasicoherent sheaves taking values in
$\R_{\ge 0}\cup\{+\infty\}$ and satisfying the following two properties:

\noindent
(i) $\rk$ is additive in exact triples;

\noindent
(ii) If a quasicoherent sheaf is represented as a filtering
union of quasicoherent subsheaves, $M=\cup_{i\in I}M_i$, then
$$\rk M=\lim_{i\in I} \rk M_i.$$
\end{thm}

Actually, the construction works in the following general framework (we keep
the notations and conventions of section \ref{ind-sec}).
Let $\AA$ be an abelian category and let $\Ind\AA$ be the corresponding
category of ind-objects of $\AA$. 
Let $\rk:K_0(\AA)\to\R$ be a homomorphism such that
$\rk(A)>0$ for any nonzero object $A\in\AA$.

\begin{thm}\label{gen-rank-thm}
There exists a unique extension of the function
$\rk$ from objects of $\AA$ to objects of $\Ind\AA$ taking values in
$\R_{\ge 0}\cup\{+\infty\}$ and satisfying the following two properties:

\noindent
(i) $\rk$ is additive in exact triples;

\noindent
(ii) If an ind-object $X$ is represented as a union of subobjects, $X=\cup_{i\in I}X_i$, then
$$\rk X=\lim_{i\in I} \rk X_i.$$
\end{thm} 

Note that Theorem \ref{rank-thm} is an immediate 
consequence of this theorem.

First, we are going to define the rank of $\AA$-generated ind-objects (see section \ref{ind-sec}).
By Proposition \ref{fingen-adm-prop} such an ind-object can be presented in the form
$P/S$, where $P\in\AA$, and $S$ is a union of subobjects of $P$ in $\AA$.
Therefore, it is natural to make the following definition.

\begin{defi} Let $S\in\Ind\AA$ be an ind-subobject of $P\in\AA$. 
Then we set 
$$\rk_P(S)=\sup\{\rk S'\ |\ S'\sub S, S'\in\AA\},$$
i.e., we define the rank of $S$ in $P$ as the supremum of the ranks of
all subobjects of $P$ in $\AA$ contained in $S$.
\end{defi}

Note that if $S$ is itself in $\AA$ then $\rk_P(S)=\rk S$.
It is clear from the definition that $\rk_P(S)\le \rk(P)$ for every
$S\sub P$ and that $\rk_P$ is monotone with respect to
inclusions. This function also satisfies the following continuity condition. 

\begin{lem}\label{rk-cont-lem} 
Let $(S_i)_{i\in I}$ be a filtering collection 
of ind-subobjects of $P\in\AA$. 
Then for $S=\cup_{i\in I} S_i$ we have
$$\rk_P(S)=\lim_{i\in I}\rk_P(S_i).$$
\end{lem}

\Pf . It is clear that $\rk_P(S)\ge\sup\{\rk_P(S_i)\ |\ i\in I\}$.
On the other hand, if $S'\sub S$ is a subobject such that $S'\in\AA$,
then there exists $i\in I$ such that $S'\sub S_i$ (by 
Proposition \ref{ind-prop}(a)). This implies that
$\rk_P(S)\le \sup\{\rk_P(S_i)\ |\ i\in I\}$.
Hence, 
$$\rk_P(S)=\sup\{\rk_P(S_i)\ |\ i\in I\}.$$
Since the function $\rk_P$ is monotone, 
we can replace $\sup$ with $\lim$.
\ed

The above lemma also implies that in the definition of $\rk_P(S)$
it suffices to take the supremum over
any collection of subobjects of $P$ in $\AA$ whose union is $S$.

\begin{prop}\label{fingen-rk-prop} 
Let $X=P/S$ be an $\AA$-generated ind-object, where $P\in\AA$. 
Then the nonnegative real number
$$\rk X:=\rk P-\rk_P(S)$$ 
does not depend on a presentation of $X$ in the form $P/S$.
\end{prop}

\Pf . If $P\to X$ and $P'\to X$ are surjections (where $P,P'\in\AA$) 
then they are dominated
by the surjection $P\oplus P'\to X$. Hence, it suffices to compare 
presentations $P/S$ and $P'/S'$ of $X$ in the case $P'\sub P$.
In this situation $S'=S\cap P'$ and $P=S+P'$. 
By assumption we have $S=\cup_{i\in I} S_i$, where $(S_i)$ is a filtering
collection of subobjects of $S$ contained in $\AA$.
Then $S'=\cup_{i\in I} S_i\cap P'$, so
$$\rk_P(S)=\sup\{\rk S_i\ |\ i\in I\},$$
$$\rk_{P'}(S')=\sup\{\rk S_i\cap P'\ |\ i\in I\}.$$
Now we observe that since $P\in\AA$, there exists $i_0\in I$
such that $P=S_{i_0}+P'$. It follows that for all $S_i\supset S_{i_0}$ we
have $\rk S_i\cap P'=\rk S_i+\rk P'-\rk P$. Hence,
$\rk_{P'}(S')=\rk_P(S)+\rk P'-\rk P$.
\ed

At this point we will only check the following expected property
of the rank function on $\AA$-generated ind-objects.

\begin{lem}\label{monotone-lem} 
Let $X'\sub X$ be an embedding of $\AA$-generated 
ind-objects. Then $\rk X'\le \rk X$.
\end{lem}

\Pf . We can find surjections $P\to X$ and $P'\to X'$ such that
$P,P'\in\AA$ and $P'\sub P$.
Let $S$ and $S'$ be kernels of these maps, so that $X=P/S$, $X'=P'/S'$.
Then $S'=S\cap P'$. Let $(S_i)_{i\in I}$ be a filtering collection
of subobjects of $S$ contained in $\AA$, such that $S=\cup_{i\in I}S_i$. 
As in the proof of Proposition
\ref{fingen-rk-prop} we see that
$$\rk_{P'}(S')=\sup\{\rk S_i\cap P'\ |\ i\in I\}.$$
Since for every $i\in I$ we have an embedding
$S_i/S_i\cap P'\sub P/P'$, it follows that
$$\rk S_i-\rk S_i\cap P'\le \rk P-\rk P'.$$
Passing to the limit in $i\in I$ we derive
that 
$$\rk_P(S)-\rk_{P'}(S')\le\rk P-\rk P'$$ 
which is equivalent to the desired inequality.
\ed

Since every ind-object $X$ is the union of a filtering collection
of its finitely generated ind-subobjects, we can now define the rank of an arbitrary ind-object
$X$ by setting
$$\rk X=\sup\{\rk X^f\ |\ X^f\sub X,\ X^f\text{ is $\AA$-generated}\}.$$
By Lemma \ref{monotone-lem}, if $X$ is itself $\AA$-generated then
this definition agrees with the old one.
We also have the following analogue of Lemma \ref{rk-cont-lem}.

\begin{lem}\label{rk-cont-lem2} 
Let $(X_i)$ be a filtering collection of
subobjects in $X\in\Ind\AA$ such that $X=\cup_{i\in I}X_i$. Then
$$\rk X=\lim_{i\in I}\rk X_i.$$
\end{lem}

\Pf . The proof follows the proof of Lemma \ref{rk-cont-lem}
step by step (recall that by Lemma \ref{monotone-lem}
the rank function on $\AA$-generated ind-objects is monotone).
\ed

We need one more lemma for the proof of Theorem \ref{gen-rank-thm}.

\begin{lem}\label{fingen-triple-lem} 
Let $0\to A\to B\to C\to 0$ be an exact triple of
ind-objects such that $B$ is $\AA$-generated.
Then $\rk B=\rk A+\rk C$.
\end{lem} 

\Pf . Let $B=P/S$, $C=P/T$, where $P$ is an object of $\AA$ and
$S\sub T\sub P$ are its ind-subobjects. We have $\rk B=\rk P-\rk_P(S)$,
$\rk C=\rk P-\rk_P(T)$. Let $(T_i)_{i\in I}$ be a filtering collection
of subobjects in $T$ contained in $\AA$, such that $T=\cup_{i\in I}T_i$.
Then $A=T/S=\cup_{i\in I}T_i/T_i\cap S$, so by Lemma \ref{rk-cont-lem2}
we obtain
$$\rk A=\lim_{i\in I}\rk T_i/T_i\cap S.$$ 
By definition, we have 
\begin{equation}\label{rkTi-eq}
\rk T_i/T_i\cap S=\rk T_i-\rk_{T_i}(T_i\cap S)=\rk T_i-\rk_P(T_i\cap S).
\end{equation}
Since $S=\cup_{i\in I}T_i\cap S$, by Lemma \ref{rk-cont-lem}
we get 
$$\lim_{i\in I}\rk_P(T_i\cap S)=\rk_P(S).$$
Therefore, passing to the limit in \eqref{rkTi-eq} we obtain
$$\rk A=\rk_P(T)-\rk_P(S).$$
\ed

\noindent
{\it Proof of Theorem \ref{gen-rank-thm}.}
It is easy to see that any extension of $\rk$ satisfying (i) and (ii)
should coincide with the rank function constructed above. Note also
that our rank function satisfies (ii) by Lemma \ref{rk-cont-lem2}. 
It remains to check its additivity in exact triples. 
Let $0\to A\to B\to C\to 0$
be an exact triple in $\Ind\AA$.
We have $B=\cup_{i\in I} B_i$, where $B_i$'s are $\AA$-generated subobjects.
Let $C_i$ be the image of $B_i$ under the map to $C$ and let
$A_i=A\cap B_i$. Then $C=\cup_{i\in I}C_i$ and $A=\cup_{i\in I}A_i$.
By Lemma \ref{fingen-triple-lem} for each $i$ we have
$$\rk B_i=\rk A_i+\rk C_i.$$
Passing to the limit and using Lemma \ref{rk-cont-lem2} we derive
that $\rk B=\rk A+\rk C$.
\ed

\subsection{Quasi vector bundles}\label{quasi-bun-sec}

Recall that these are quasicoherent sheaves that are filtering unions of holomorphic bundles.

\begin{lem}\label{fingen-quasi-lem} 
Every finitely generated holomorphic submodule of a quasi vector bundle
is a vector bundle. 
\end{lem}

\Pf . Indeed, let $M=\cup_i M_i$, where $M_i$ are vector bundles.
Then every finitely generated submodule of $M$ is contained in some $M_i$,
hence itself is a vector bundle by Corollary \ref{fingen-proj-cor}.
\ed

\begin{lem}\label{quasi-sub-lem} 
Let $M$ be a quasi vector bundle, $N\sub M$ a quasicoherent subsheaf.
Then $N$ is a quasi vector bundle.
\end{lem}

\Pf . Since $N$ is a union of its finitely generated holomorphic submodules,
this follows immediately from the previous lemma. 
\ed

\begin{prop}\label{quasi-res-prop} 
Every quasicoherent sheaf on $T$ can be represented in the form
$P_1/P_0$, where $P_0\sub P_1$ are quasi vector bundles.
\end{prop}

\Pf . Recall that by Lemma \ref{quasi-sub-lem}
a quasicoherent subsheaf of a quasi vector bundle is itself a quasi vector bundle.
Therefore, it is enough to prove
that for every quasicoherent sheaf $M$ there exists a quasi vector
bundle $P$ and a surjection $P\to M$. By Proposition \ref{adm-prop2}
we have $M=\cup_{i\in I} M_i$ with $M_i=P_i/Q_i$, where $P_i$ are vector bundles.
Set $P=\oplus_{i\in I} P_i$ and define the morphism $P\to M$ using the natural
morphisms $P_i\to M_i\sub M$. It is clear that this morphism is surjective and
that $P$ is a quasi vector bundle. 
\ed

For later use we record here 
one more simple observation.

\begin{lem}\label{quot-quasi-lem}
Let $V\sub M$ be an embedding of a vector bundle into a quasi vector bundle.
Then $M/V$ is a quasi vector bundle.
\end{lem}

\Pf . Note that if $M'\sub M$ is a vector subbundle then $V+M'\sub M$ is
still a vector bundle, e.g., by Lemma \ref{fingen-quasi-lem}.
It follows that $M$ can be represented
as a union of vector bundles $M=\cup_i M_i$, where $V\sub M_i$ for all $i$. Hence,
$M/V=\cup_i M_i/V$.
\ed

\subsection{Torsion and torsion free sheaves}
The following two classes of quasicoherent sheaves will be also important for us.

\begin{defi}
Let $M$ be a quasicoherent sheaf on $T$. We say that $M$ is a 
{\it torsion sheaf} if $\rk M=0$. We say that $M$ is a
{\it torsion-free sheaf} if for every nonzero submodule $M'\sub M$
one has $\rk M'>0$.
\end{defi}

\begin{prop} (i) Let $M$ be a torsion sheaf and $N$ be a torsion-free sheaf.
Then $\Hom(M,N)=0$.

\noindent (ii) For every quasicoherent sheaf $M$ there exists maximal
torsion subsheaf $M_{tors}\sub M$. The quotient $M/M_{tors}$ is
torsion free.
\end{prop}

\Pf . (i) If $f:M\to N$ is a morphism then $\im(f)$ has rank $0$. But
it is a subsheaf of $N$, hence, $\im(f)=0$.

\noindent
(ii) If $N_1$ and $N_2$ are torsion subsheaves in $M$. Then there exists
a surjection $N_1\oplus N_2\to N_1+N_2\sub M$. Hence, $\rk(N_1+N_2)=0$,
so $N_1+N_2$ is also a torsion sheaf. It follows that the union of
all torsion subsheaves in $M$ is itself a subsheaf $M_{tors}\sub M$.
Furthermore, by Lemma \ref{rk-cont-lem2} we have $\rk M_{tors}=0$.
If $N\sub M/M_{tors}$ is torsion subsheaf then its preimage in
$M$ is also a torsion sheaf by additivity of the rank. Hence, $N=0$.
\ed

\begin{prop} Every quasi vector bundle is torsion free.
\end{prop}

\Pf . Let $V$ be a quasi vector bundle. It suffices to prove that
for every finitely generated quasicoherent subsheaf $W\sub V$ with
$\rk W=0$ one has $W=0$. But this follows immediately from 
Lemma \ref{fingen-quasi-lem}.
\ed

Modifying slightly the proof of Theorem \ref{main-thm} we get the following result.

\begin{thm}\label{quasi-main-thm} 
For every stable bundle $P$ and
every real number $r$ such that $0<r<\rk P$ there exists a countably generated
quasicoherent
subsheaf $Q\sub P$ such that $\rk Q=r$.
\end{thm}

\Pf . Using Morita equivalences we reduce to the case $\rk P=1$.
Then we can apply the construction of Theorem \ref{main-thm}
to construct a family of subbundles $V_{0,a}\sub P$ for
$a\in\BB_{\th}$, where $\rk V_{0,a}=a$ and $V_{0,a}\sub V_{0,a'}$ for
$a<a'$. Let $(a_n)$ be an increasing sequence of numbers in 
$\BB_{\th}$ such that $\lim_{n\to\infty} a_n=r$. Then we can take
$Q=\cup_{n} V_{0,a_n}$.
\ed 

Later we will see that in fact all quasicoherent subsheaves of $P$ are countably generated
(see Theorem \ref{count-thm}(i)).

\begin{cor}\label{quasi-main-cor} 
For every real number $r>0$ 
there exist a finitely generated torsion free quasicoherent sheaf $M$ of rank $r$ which is
not a quasi vector bundle.
\end{cor}

\Pf . Let $P$ be a stable bundle with $\rk P>r$. In the case
$r\in\Z\th+\Z$ we also require that $\chi(\rk P,r)<0$. 
Then by the above theorem
there exists a quasicoherent subsheaf $S\sub P$ with $\rk S=\rk P-r$. Now we define
$M$ to be the quotient of $P/S$ by its torsion part $(P/S)_{tors}$.
Then $M$ is a finitely generated torsion free sheaf of rank $r$. We claim
that $M$ is not a vector bundle. Indeed, if $r\not\in\Z\th+\Z$
then this is clear. Otherwise, using the fact that
$M$ is a quotient of a stable bundle $P$ we get a contradiction with the condition
$\chi(\rk P,r)<0$. Finally, Lemma \ref{fingen-quasi-lem} implies that $M$ is not a quasi vector bundle.
\ed

\subsection{Harder-Narasimhan filtration for quasicoherent subsheaves of vector bundles}
\label{HN-sec}

In this section we will show that every quasicoherent subsheaf of a vector bundle on $T$
has a canonical exhaustive filtration similar to the Harder-Narasimhan filtration of vector bundles. 

\begin{lem} For $N>0$ consider the subset $\MM_N\sub\R$ defined by
$$\MM_N=\{\frac{m}{m\th+n}\ |\ (m,n)\in\Z, m\le 0, 0<m\th+n<N\}.$$
Then for every $c>0$ the set $\MM_{\th}\cap [-c,0]$ is finite.
In particular, every nonempty subset of $\MM_{\th}$ has a maximal element.
\end{lem}

\Pf . This follows immediately from the fact that there are at most $N$ numbers of the form
$m\th+n$ in the interval $(0,N)$ with a given value of $m$ and that
for such a number we have $|m|/N\le |m|/(m\th+n)$.
\ed

\begin{lem} Let $P$ be a vector bundle and $\FF\sub P$ a quasicoherent
subsheaf. Then among subbundles $Q\sub P$ such that $Q\sub\FF$ there
exists a unique maximal element (by inclusion) of maximal slope. Moreover,
such $Q$ is semistable.
\end{lem}

\Pf . It is easy to see that $P$ can be embedded into a vector bundle of the form
$P_0^{\oplus N}$, where $P_0$ is a stable vector bundle, so we can replace $P$
with $P_0^{\oplus N}$. Furthermore, applying Morita equivalence we can assume
that $\rk P_0=1$. In this case the slope of any subbundle of $P_0^{\oplus N}$ belongs
to the subset $\MM_N$ considered in the above lemma. Hence, 
there exists a subbundle $Q\sub \FF$ of a maximal slope $\mu_{max}$. 
Note that such $Q$ is automatically
semisimple. Therefore, among all subbundles of maximal slope there exists a maximal
element by inclusion. It remains to prove that it is unique. But if $Q'\sub \FF$ is another
subbundle with this property then the exact sequence
$$0\to Q\to Q+Q'\to Q'/Q\cap Q'\to 0$$
shows that the $\mu(Q+Q')\ge\mu_{max}$. Indeed, it suffices to see that 
$\mu(Q'/Q\cap Q')\ge\mu_{max}$ but this follows from the semistability of $Q'$.
\ed

\begin{prop}\label{HN-quasi-prop} 
Let $P$ be a vector bundle and $\FF\sub P$ a quasicoherent subsheaf.
Then there exists a sequence of subbundles $0=F_0\sub F_1\sub F_2\sub\ldots\sub P$
such that $\FF=\cup_n F_n$ and $F_n/F_{n-1}$ is the maximal subbundle
of maximal slope in $\FF/F_{n-1}\sub P/F_{n-1}$.
\end{prop}

\Pf . The previous lemma allows to define the sequence of subbundles
$0=F_0\sub F_1\sub F_2\sub\ldots\sub P$ such that $F_n/F_{n-1}$ is the maximal
subbundle of maximal slope in $\FF/F_{n-1}\sub P/F_{n-1}$. It remains to show
that $\FF=\cup_n F_n$. If $\FF$ is a subbundle then
$(F_n)$ coincides with its usual Harder-Narasimhan filtration (as an object in the derived category 
$D^b(E)$), so in this case our claim holds. Thus, we can assume that $\FF$ is not coherent, so
that $F_n/F_{n-1}\neq 0$ for all $n>0$. 
Note that $\rk(F_n/F_{n-1})\to 0$ as $n\to\infty$. Also, we have $\mu(F_1/F_0)>\mu(F_2/F_1)>\ldots$.
We claim that this implies that $\mu(F_n/F_{n-1})\to-\infty$ as $n\to\infty$. Indeed, if
the sequence $(\mu(F_n/F_{n-1}))$ were bounded then we would have
$\deg(F_n/F_{n-1})\to 0$ as $n\to\infty$. But $\deg(F_n/F_{n-1})$ can take value $0$ only once,
so we get a contradiction that proves our claim.
Now let $Q\sub P$ be an arbitrary subbundle such that $Q\sub\FF$.
We have to prove that $Q\sub F_n$ for $n\gg 0$. 
Since $Q$ is a direct sum of semistable subbundles, it suffices to consider the case when
$Q$ is semistable. Choose $n$ such that $\mu(Q)>\mu(F_n/F_{n-1})$.
Assume that $Q\not\sub F_{n-1}$. Then
$Q'=Q/(Q\cap F_{n-1})$ is a subbundle of $P/F_{n-1}$ contained in $\FF/F_{n-1}$.
Furthermore, since $Q$ is semistable, we have
$\mu(Q')\ge\mu(Q)>\mu(F_n/F_{n-1})$. But this contradicts to the assumption that
$F_n/F_{n-1}$ has maximal slope among all subbundles of $P/F_{n-1}$ containing $\FF/F_{n-1}$.
Therefore, $Q$ is contained in $F_{n-1}$.
\ed

\subsection{Countably generated quasicoherent sheaves}
\label{count-sec}

We say that a quasicoherent sheaf $M$ is {\it countably generated} if its underlying
$A_{\th}$-module has a countable set of generators. 
Equivalently, $M=\cup_{n\ge 1}M_n$ for a chain $M_1\sub M_2\sub\ldots\sub M$
of finitely generated quasicoherent subsheaves. Note that for such $M$
there exists a surjection $P\to M$, where $P$ is a countably generated quasi vector bundle
(pick surjections $P_i\to M_i$ and set $P=\oplus P_i$).
The importance of countability is due to the following general result.

\begin{lem}\label{proj-obj-lem} 
Let $P$ be an object of an abelian category $\AA$ such that
$P=\dlim_i P_i$ for some inductive system $(P_i)$ with a countable set of indices,
such that all
objects $P_i$ are projective and every arrow $P_i\to P_j$ is an
embedding of a direct summand. Then $P$ itself is projective.
\end{lem}

\Pf . By assumption, for every arrow $P_i\to P_j$ and
every $A\in\AA$ the morphism
$$\Hom_{\AA}(P_j,A)\to\Hom_{\AA}(P_i,A)$$
is surjective. Therefore, the functor
$$A\to\Hom_{\AA}(P,A)=\projlim_i \Hom_{\AA}(P_i,A)$$
is exact (since the Mittag-Leffler condition is satisfied, see \cite{EGA}, ch.~0, \S~13), so
$P$ is projective.
\ed

\begin{prop}\label{quasi-proj-mod-prop} 
A countably generated quasi vector bundle is projective as an $A_{\th}$-module.
\end{prop}

\Pf . This follows immediately from Lemma \ref{proj-obj-lem} since every
embedding of holomorphic vector bundles is an embedding of a direct summand
on the level of $A_{\th}$-modules.
\ed

\begin{lem}\label{count-lem} 
Let $M$ be a countably generated quasicoherent sheaf,
$N\sub M$ be a quasicoherent subsheaf. Then $N$ is countably generated.
\end{lem}

\Pf . Consider a surjection $P\to M$, where $P$ is a countably generated quasi vector bundle.
It is enough to prove that every quasicoherent subsheaf $N\sub P$ is countably generated.
Furthermore, if $P$ is a union of a sequence of vector bundles $P_1\sub P_2\sub\ldots $
then $N=\cup_i (N\cap P_i)$. Therefore, it suffices to consider the case when $P$ is a vector bundle.
Then the assertion follows from Proposition \ref{HN-quasi-prop}.
\ed

\begin{thm}\label{count-thm} 
(i) The full subcategory $\Qcoh^c(T)\sub\Qcoh(T)$ of countably generated sheaves
is closed under passing to subobjects and quotients, and under extensions.
In other words, $\Qcoh^c(T)$ is a Serre subcategory in $\Qcoh(T)$.

\noindent
(ii) For any $M\in\Qcoh^c(T)$ the projective dimension
of the underlying $A_{\th}$-module is at most $1$.

\noindent
(iii) A sheaf $M\in\Qcoh^c(T)$ is a quasi vector bundle
iff the underlying $A_{\th}$-module is projective.
\end{thm}

\Pf . (i) It is clear that $\Qcoh^c(T)$ is closed under passing to quotients. The assertion about
subobjects is Lemma \ref{count-lem}. The assertion about extensions is clear.

\noindent
(ii) This follows from Proposition \ref{quasi-proj-mod-prop}.

\noindent
(iii) The ``only if" part is Proposition \ref{quasi-proj-mod-prop}.
Conversely, assume that $M\in\Qcoh^c(T)$
is projective as an $A_{\th}$-module. 
Let $M'\sub M$ be a finitely generated quasicoherent subsheaf. Then
$M/M'$ is an object of $\Qcoh^c(T)$. Hence, by part (ii)
$M'$ is still projective as an $A_{\th}$-module. Therefore, $M'$ is
a holomorphic vector bundle. Since $M$ is a union of finitely generated subsheaves,
it is a quasi vector bundle.
\ed

\begin{cor}\label{count-ideal-cor} 
Any quasicoherent ideal $I\sub A_{\th}$ is countably generated.
\end{cor}

\section{Sheaves at the general point of a noncommutative torus}
\label{quasi-gen-sec}

\subsection{Quasi vector bundles at the general point}

Let $\Tors\sub\Qcoh(T)$ be the full subcategory consisting of
torsion sheaves $M$ (i.e., sheaves with $\rk M=0$). This is a Serre subcategory
of $\Qcoh(T)$, so we can consider the quotient-category
$$\Qcoh(\eta_T)=\Qcoh(T)/\Tors$$
which is a noncommutative 
analogue of the category of quasicoherent sheaves on a general
point of an elliptic curve. Note that $\Qcoh(\eta_T)$
is a $\C$-linear abelian category and there is a canonical exact functor
$\Qcoh(T)\to\Qcoh(\eta_T)$.

\begin{prop} Let $P_1$ and $P_2$ be quasicoherent sheaves on $T$. 
Assume that $P_2$ is torsion free. Then the natural
morphism
$$\Hom_{\Qcoh(T)}(P_1,P_2)\to\Hom_{\Qcoh(\eta_T)}(P_1,P_2)$$
is injective.
\end{prop}

\Pf . 
By definition, a morphism from $P_1$ to $P_2$ in $\Qcoh(\eta_T)=\Qcoh(T)/\Tors$
is given by a morphism $P'_1\to P_2/F$, where $P'_1\sub P$ and $F\sub P_2$
are quasicoherent subsheaves such that $\rk P_1/P'_1=\rk F=0$. But $P_2$ is torsion free, 
hence $F=0$. Thus, 
$$\Hom_{\Qcoh(\eta_T)}(P_1,P_2)=\dlim_{P'_1\sub P_1: \rk P_1/P'_1=0}
\Hom_{\Qcoh(T)}(P'_1,P_2).$$
It remains to check that if a morphism $f:P_1\to P_2$ vanishes on a subsheaf
$P'_1\sub P_1$ such that $\rk P_1/P'_1=0$ then $f=0$. But such $f$ factors
through a morphism $P_1/P'_1\to P_2$. Since $\rk P_1/P'_1=0$ and $P_2$ is 
torsion free, such a morphism has to be zero.
\ed

\begin{cor} The functor $\Vect(T)\to\Qcoh(\eta_T)$ is faithful.
\end{cor}

\begin{lem}\label{surj-lim-lem}
Let $f:M\to M'$ be a surjection in $\Qcoh(T)$. Then there exist inductive systems
$(M_i)$ and $(M'_i)$ in $\Vect(T)$ such that $\dlim M_i\simeq M$, $\dlim M'_i\simeq 
M'$,
and a morphism of inductive systems $(M_i)\to (M'_i)$ inducing $f$, such that
every morphism $M_i\to M'_i$ is a surjection.
\end{lem}

\Pf . By Proposition \ref{quasi-res-prop} we can find a quasi vector bundle $P$ and a surjection $P\to M$.
Thus, we can assume that $f$ is a natural morphism $P/S\to P/S'$, where $S\sub S'\sub P$
are subsheaves. 
Let $S=\cup_{i\in I} S_i$ (resp., $S'=\cup_{j\in J} S'_j$), where
$S_i\sub P$ (resp., $S'_j\sub P$) are holomorphic vector bundles. We can assume that the
sets of indices $I$ and $J$ are the same (e.g., replacing both by $I\times J$).
Furthermore, replacing $S'_i$ with $S'_i+S_i$
(which is still a subbundle by Lemma \ref{fingen-quasi-lem}) we can assume
that $S_i\sub S'_i$. Then $(P/S_i)\to (P/S'_i)$ is the required morphism of inductive systems.
\ed

\begin{lem}\label{lifting-lem} 
Assume that we have a commutative diagram in $\Qcoh(T)$ of the form
\begin{diagram}
M &\rTo^{f}&M'\\
\uTo{g}&&\uTo{}\\
S &\rTo^{i} &P
\end{diagram}
where $S$ is a vector bundle, $P$ is a quasi vector bundle of finite rank,
$f$ is a surjection, $i$ is an embedding. Then for every $\eps>0$ there
exists a vector bundle $Q\sub P$ such that $S\sub Q$, $\rk Q>\rk P-\eps$
and there exists a morphism $Q\to M$ making the following diagram commutative:
\begin{diagram}
&& M &\rTo^{f}&M'\\
&\ruTo{g}&\uTo{}&&\uTo{}\\
S &\rTo{}&Q&\rTo{}&P
\end{diagram}
\end{lem}

\Pf . We split the proof in two steps.

\noindent 
{\bf Step 1}. {\it Assume that $M$ and $M'$ are vector bundles}.
Also, without loss of generality we can assume also that $P$ is a vector bundle.
Indeed, otherwise we can replace $\eps$ with $\eps/2$
and $P$ with some bundle $P'\sub P$ such that $S\sub P'$ and $\rk P'>\rk P-
\eps/2$.
Furthermore, replacing $M$ with the fibered product of $M$ and $P$ over $M'$
we can assume that $P=M'$. Let $N=\ker(f)$. Then we have an exact
sequence
$$0\to N\to M/S\to P/S\to 0.$$
Using Lemma \ref{eps-lem} we can find
a subbundle $Q'\sub P/S$ such that $\rk Q'>\rk P/S-\eps$ and 
$\Ext^1(Q',N)=0$. Then the pull-back of the above exact sequence
to $Q'\sub P/S$ splits. Let $Q'\to M/S$ be a splitting and let
$Q\sub P$ be the preimage of $Q'$ in $P$. Since $M$ is the
fibered product of $M/S$ and $P$ over $P/S$ we obtain a morphism
$Q\to M$ with required properties. 

\noindent 
{\bf Step 2}. Now 
using Lemma \ref{surj-lim-lem} we can find inductive systems of
vector bundles $(M_i)$ and $(M'_i)$ and a system of surjections
$f_i:M_i\to M'_i$ inducing $f$. Since the functor
$\Hom(S,-)$ on $\Qcoh(T)$ commutes with inductive limits, there exists
a commutative diagram of the form
\begin{diagram}
M_i &\rTo^{f_i}&M'_i\\
\uTo{}&&\uTo{}\\
S &\rTo^{i} &P
\end{diagram}
inducing our original diagram. It remains to apply Step 1.
\ed

\begin{lem}\label{count-sur-lem} 
Let $M$ be a countably generated quasicoherent sheaf of finite rank.
Then there exists a quasi vector bundle $P$ of finite rank and a surjection
$P\to M$.
\end{lem} 

\Pf . Since $M$ is countably generated, there exists a sequence of finitely
generated subsheaves $M_1\sub M_2\sub\ldots\sub M$ such that $M=\cup_{n\ge 1}M_n$.
Furthermore, we can choose a sequence of vector bundles $P_1\sub P_2\sub\ldots$
and of compatible surjections $f_n:P_n\to M_n$. Let $K_n=\ker f_n$.
We are going to choose recursively a sequence of vector bundles
$Q_n\sub K_n$ such that $Q_n=Q_{n+1}\cap P_n$ and 
$\rk K_n-\rk Q_n<n/(n+1)$. For $n=1$ we choose $Q_1$ to be any subbundle
of $K_1$ such that $\rk K_1-\rk Q_1<1/2$.
Assume that $Q_n$ is already constructed and let us set $P'_n=P_n/Q_n$,
$P'_{n+1}=P_{n+1}/Q_n$, $K'_n=K_n/Q_n\sub P'_n$ and 
$K'_{n+1}=K_{n+1}/Q_n\sub P'_{n+1}$. Note that $K'_{n+1}\cap P'_n=K'_n$
and $\rk K'_n<n/(n+1)$.
For every $\eps>0$ 
we can choose a vector bundle $R\sub K'_{n+1}\sub P'_{n+1}$ such that 
$\rk R>\rk K'_{n+1}-\eps$. Applying Lemma \ref{lifting-lem}
to the surjection $R\to R/R\cap P'_n$ we find a subbundle
$Q'_{n+1}\sub R\cap P'_n$ such that $Q'_{n+1}$ lifts to a subbundle of $R$ and
$\rk Q'_{n+1}>\rk R/R\cap P'_n-\eps$. Note that
$R\cap P'_n\sub K'_n$. Hence, $\rk R/R\cap P'_n>\rk R-n/(n+1)$, and therefore,
$$\rk Q'_{n+1}>\rk R-n/(n+1)-\eps.$$
Viewing $Q'_{n+1}$ as a subbundle
$Q'_{n+1}\sub R\sub K'_{n+1}$ 
let us define $Q_{n+1}$ as a preimage of $Q'_{n+1}$ in $K_{n+1}$. 
Since $Q'_{n+1}\cap P'_n=0$, we obtain $Q_{n+1}\cap P_n=Q_n$.
Also, 
$$\rk K_{n+1}-\rk Q_{n+1}=\rk K'_{n+1}-\rk Q'_{n+1}<\rk R+\eps-\rk Q'_{n+1}<
n/(n+1)+2\eps.$$
Thus, if we choose $\eps$ sufficiently small we will satisfy the condition
$\rk K_{n+1}-\rk Q_{n+1}<(n+1)/(n+2)$.

Now let us consider the sequence of vector bundles 
$\ov{P}_1\sub\ov{P}_2\sub\ldots$,
where $\ov{P}_n=P_n/Q_n$. Set $\ov{P}=\cup_{n\ge 1}\ov{P}_n$. By definition
$\ov{P}$ is a quasi vector bundle. Furthermore, 
$$\rk\ov{P}_n=\rk P_n-\rk Q_n<\rk P_n-\rk K_n+n/(n+1)=\rk M_n+n/(n+1).$$
Hence, $\rk\ov{P}\le \rk M+1$.
\ed

\begin{thm}\label{proj-thm}
(i) A quasicoherent sheaf of finite rank is a projective object of $\Qcoh(\eta_T)$  
(resp., $\Qcoh^f(\eta_T)$) iff it is isomorphic
to a quasi vector bundle.

\noindent
(ii) The categories $\Qcoh(\eta_T)$ and $\Qcoh^f(\eta_T)$ have enough projective objects.
The cohomological dimension of $\Qcoh^f(\eta_T)$ is at most $1$.
\end{thm}

\Pf . (i) First, let us prove that a vector bundle $P$ considered as
an object of $\Qcoh(\eta_T)$ is projective.
Every surjection in $\Qcoh(\eta_T)$ can be represented by a morphism
$f:M\to M'$ in $\Qcoh(T)$ such that $\rk \coker(f)=0$. We have to show
that every morphism from $P$ to $M'$ in $\Qcoh(\eta_T)$ factors through $f$.
By definition, every such morphism is given by a morphism
$P'\to M'/F$, where $P'\sub P$ and $F\sub M'$ are such that $\rk P/P'=0$ and
$\rk F=0$. Replacing $M'$ by $M'/F$ we can assume that $F=0$. Also by
Lemma \ref{quasi-sub-lem} we can replace $P$ by $P'$. Thus, it suffices to prove
that every morphism $P\to M'$ in $\Qcoh(T)$ factors through $f$ in 
$\Qcoh(\eta_T)$.
Let $P'\sub P$ be the preimage of $\im(f)\sub M'$. Then $\rk P/P'=0$, so 
replacing $M'$ by $\im(f)$ and $P$ by $P'$ we can assume that $f$ is surjective.
Iterating Lemma \ref{lifting-lem} we can construct a sequence of bundles
$S_1\sub S_2\sub\ldots\sub P$ such that $\rk S_n>\rk P-1/n$ equipped with 
a system of compatible liftings of the induced morphisms $S_n\to M'$ to morphisms
$S_n\to M$. Indeed, to construct $S_1$ we apply Lemma \ref{lifting-lem} with $S=0$ 
and
$\eps=1$, and
set $S_1=Q$. If $S_n$ is already constructed then we apply Lemma \ref{lifting-lem}
with $S=S_n$ and $\eps=1/(n+1)$, and set $S_{n+1}=Q$.
Let $P'=\cup_n S_n$. Then the induced morphism $P'\to M'$ factors through $f$
and $\rk P/P'=0$, so we are done.

Now if $P$ is any quasi vector bundle of finite rank
then we can choose a sequence of vector bundles $P_1\sub P_2\sub\ldots\sub P$
such that $\lim \rk P_n=\rk P$. Note that $\cup P_n\simeq P$ in $\Qcoh(\eta_T)$,
so we can assume that $P=\cup _n$.
As we have seen above, every $P_i$ is a projective object in $\Qcoh(\eta_T)$.
Therefore, by Lemma \ref{proj-obj-lem} $P$ is also projective.

Conversely, assume that $M$ is a quasicoherent sheaf of finite rank which is a projective
object of $\Qcoh^f(\eta_T)$. Without loss of generality 
we can assume that $M$ has no torsion.
Also, we can choose a sequence of finitely generated
subsheaves $M_1\sub M_2\sub\ldots\sub M$ such that $\lim\rk M_n=\rk M$.
Replacing $M$ by $\cup_{n\ge 1}M_n$ we can assume that $M$ is countable
generated. Then by Lemma \ref{count-sur-lem} we can find
a surjection $P\to M$ in $\Qcoh(T)$,
where $P$ is a quasi vector bundle of finite rank.
Now our assumption implies that there exists a splitting $M\to P$ in $\Qcoh(\eta_T)$.
Since $P$ has no torsion, this splitting is given by a morphism $f:M'\to P$ in
$\Qcoh(T)$, where $M'\sub M$ is such that $\rk M/M'=0$. Shrinking $M'$ if
necessary we can assume that the composition $\pi\circ f:M'\to M$
coincides with the embedding of $M'$ into $M$. Hence, $f$ is an embedding.
By Lemma \ref{quasi-sub-lem} this implies that $M'$ is a quasi vector bundle.
It remains to observe that $M'\simeq M$ in $\Qcoh(\eta_T)$. 

\noindent
(ii) Let $M$ be any quasicoherent sheaf. Then we can find a collection of
vector bundles $(P_i)$ and a surjection $P=\oplus P_i\to M$. Note that each
$P_i$ is a projective object in $\Qcoh(\eta_T)$, hence $P$ is also projective.
This shows that $\Qcoh(\eta_T)$ has enough projective objects.
The similar statement for the subcategory $\Qcoh^f(\eta_T)$ follows from   
Lemma \ref{count-sur-lem}. Now the fact that the cohomological dimension of $\Qcoh^f(\eta_T)$
is $\le 1$ follows easily from Lemma \ref{quasi-sub-lem} and from part (i).
\ed

\subsection{Isomorphisms in $\Qcoh(\eta_T)$}

\begin{lem}\label{rank-bun-isom-lem} 
Let $P$ and $P'$ be holomorphic vector bundles such that $\rk P=\rk P'$.
Then $P\simeq P'$ in $\Qcoh(\eta_T)$.
\end{lem}

\Pf . By Theorem \ref{def-thm} it is enough to prove that
deformation equivalent bundles 
become isomorphic in $\Qcoh(\eta_T)$. Note that
Theorem \ref{proj-thm} implies that every exact triple of vector bundles
splits in $\Qcoh(\eta_T)$. Hence, it is enough to prove the assertion in the case when
$P_1$ and $P_2$ are stable. Using Morita equivalences we can reduce to the case
when $\rk P_1=\rk P_2=1$. Let us apply the construction of Theorem \ref{main-thm}
to get a family $(V_{a,b})$ (resp., $V'_{a,b}$)
of stable subquotients of $P$  (resp., $P'$) numbered by
the subsegments $[a,b]$ of the division process described in section \ref{div-constr-sec}.
Let also $V_{0,a}\sub P$ (resp., $V'_{0,a}\sub P'$) be the corresponding subbundles
numbered by $a\in\BB_{\th}$.
We can choose these families in such a way that $V_{a,b}\simeq V'_{a,b}$
for all $[a,b]$ appearing in the division process. Since all exact sequences of vector bundles split in $\Qcoh(\eta_T)$, it follows that $V_{0,a}\simeq V'_{0,a}$ in $\Qcoh(\eta_T)$.
Hence, we get an isomorphism 
$$P\simeq\cup_{a\in\BB_{\th}}V_{0,a}\simeq\cup_{a\in\BB_{\th}}V'_{0,a}\simeq P'$$
in $\Qcoh(\eta_T)$.
\ed

\begin{lem}\label{sub-rank-chi-lem} 
Let $V$ be a vector bundle. Then there
exists an element $v\in(\Z+\Z\th)_{>0}$ such that for every $r\in\Z+\Z\th$ such 
that $0\le r\le\rk V$ and $\chi(r,v)\ge 0$,
there exists a subbundle $W\sub V$ with $\rk W=r$.
Furthermore, if $V$ is semistable then we can take $v=\rk V$.
\end{lem}

\Pf . If $V$ is stable then the assertion follows from Theorem \ref{main-thm}. 
In the general case let $0\sub F_1V\sub F_2V\sub\ldots\sub F_nV=V$
be a filtration such that the bundles $F_iV/F_{i-1}V$
are stable and $\mu(F_1V)\ge\mu(F_2V/F_1V)\ge\ldots\ge\mu(F_nV/F_{n-1}V)$.
Set $v_i=\rk(F_iV/F_{i-1}V)$. We claim that we can take $v=v_n$.
Indeed, for every $r$ between $0$ and $\rk V$ there exists $i$, $1\le i\le n$, 
such
that $v_1+\ldots+v_{i-1}\le r\le v_1+\ldots+v_i$.
Set $r'=r-(v_1+\ldots+v_{i-1})$. Since $\chi(v_j,v_i)\le 0$ for $j<i$ we have
$$\chi(r',v_i)=\chi(r,v_i)-\sum_{j=1}^{i-1}\chi(v_j,v_i)\ge 0.$$
Hence there exists a subbundle $W'\sub F_iV/F_{i-1}V$ with $\rk W'=r'$.
It remains to take $W$ to be the preimage of $W'$ in $F_iV\sub V$.
\ed

\begin{lem}\label{pair-sub-lem} 
Let $V_1$ and $V_2$ be vector bundles. Then for every $\eps>0$ there
exist subbundles $W_1\sub V_1$ and $W_2\sub V_2$ such that 
$$\rk W_1=\rk W_2>\min(\rk V_1,\rk V_2)-\eps.$$
\end{lem}

\Pf . Indeed, let $v_1$ and $v_2$ be elements of $(\Z+\Z\th)_{>0}$ chosen
as in Lemma \ref{sub-rank-chi-lem} for $V_1$ and $V_2$, respectively.
Without loss of generality we can assume that $\chi(v_1,v_2)\ge 0$.
Then for every $r\in\Z+\Z\th$ such that $0\le r\le\min(\rk V_1,\rk V_2)$
and $\chi(r,v_1)\ge 0$, there exists subbundles $W_1\sub V_1$ and $W_2\sub V_2$
with $\rk W_1=\rk W_2=r$. Since the set of such $r$ is dense in the interval
$[0,\min(\rk V_1,\rk V_2)]$ the assertion follows.
\ed

\begin{thm}\label{quasi-thm} Let $P$ and $P'$ be quasi vector bundles
of finite ranks on $T$ such that $\rk P=\rk P'$. Then $P\simeq P'$
in $\Qcoh(\eta_T)$. 
\end{thm}

\Pf . Let $r=\rk P=\rk P'$.
First of all, by definition of the rank and by Lemma \ref{fingen-quasi-lem}, 
for every $\eps>0$ there exist
embeddings $V\sub P$ and $V'\sub P'$ in $\Qcoh(T)$, where $V$ and $V'$ are vector 
bundles
and of rank $>r-\eps$. Applying Lemma \ref{pair-sub-lem} we find subbundles
$W\sub V$ and $W'\sub V'$ such that $\rk W=\rk W'>r-2\eps$.
Since $V/W$ and $V'/W'$ are again quasi vector bundles by Lemma \ref{quot-quasi-lem},
we can apply the same procedure to $V/W$ and $V'/W'$ again, and so on.
Taking $\eps=1/2n$ at the $n$-th step, we will construct in this way a sequence 
of subbundles $0=W_0\sub W_1\sub W_2\sub\ldots\sub P$ (resp.,
$0=W'_0\sub W'_1\sub W'_2\sub\ldots\sub P'$) such that $\rk W_n=\rk W'_n>r-1/n$. 
It follows that $\rk \cup_n W_n=\rk \cup_n W'_n=r$. Hence,
$P\simeq \cup_n W_n$ and $P'\simeq \cup_n W'_n$ in $\Qcoh(\eta_T)$.
Now applying Theorem \ref{proj-thm}(i) we derive that
$\cup_n W_n\simeq \oplus_{n\ge 1} W_n/W_{n-1}$ 
(resp., $\cup_n W'_n\simeq\oplus_{n\ge 1}W'_n/W'_{n-1}$) 
in $\Qcoh(\eta_T)$. 
It remains to observe that $W_n/W_{n-1}\simeq W'_n/W'_{n-1}$ in
$\Qcoh(\eta_T)$ for all $n\ge 1$ by Lemma \ref{rank-bun-isom-lem}. 
\ed

\begin{cor}\label{proj-isom-cor} 
Projective objects in $\Qcoh^f(\eta_T)$ are determined up to an isomorphism by their rank.
In other words, for a pair of projective objects $P,P'\in\Qcoh^f(\eta_T)$ one has
$P\simeq P'$ iff $\rk P=\rk P'$.
\end{cor}

\Pf . Combine Theorem \ref{proj-thm}(i) and Theorem \ref{quasi-thm}.
\ed

\begin{cor}\label{K0-cor} 
One has $K_0(\Qcoh^f(\eta_T))\simeq\R$ and the effective cone
is exactly $\R_{>0}\sub\R$.
\end{cor}

\subsection{Equivalences with categories of modules}

Let $P$ be a quasi vector bundle of finite rank and let $R_P=\End_{\Qcoh(\eta_T)}(P)$. 
By Theorem \ref{proj-thm} the corresponding functor 
\begin{equation}\label{functor-eq}
\Ga_P:\Qcoh(\eta_T)\to \mod-R_P: M\mapsto\Hom_{\Qcoh(\eta_T)}(P,M)
\end{equation}
is exact. Below we are going to study properties of the ring $R_P$ and of the
functor $\Ga_P$.

Recall that a ring $R$ is called {\it right semihereditary} if every finitely 
generated right ideal in $R$ is projective.

\begin{prop}\label{semiher-prop} For every quasi vector bundle $P$ of finite rank the ring $R_P$
is right semihereditary.
\end{prop}

\Pf .  A finitely generated ideal $I\sub R_P$ is the image of a morphism
of $R_P$-modules $R_P^{\oplus n}\to R_P$. Such a morphism is the image under 
$\Ga_P$
of a morphism $f:P^{\oplus n}\to P$ in $\Qcoh(\eta_T)$. By Theorem 
\ref{proj-thm} the projective dimension of $\coker(f)$ is $\le 1$, hence,
$\im(f)$ is projective. It follows that $\im(f)$ is a direct summand of
$P^{\oplus n}$. Therefore, $I\simeq\Ga_P(\im(f))$ is a direct summand of
$R_P^{\oplus n}$.
\ed

\begin{lem}\label{equiv-cat-lem} 
Let $P$ be a projective object in an abelian category $\AA$, and
let $\lan P\ran\sub\AA$
denote the full subcategory consisting of objects
that can be presented as the cokernel of a morphism of the form
$P^{\oplus m}\to P^{\oplus n}$. Then the functor $\Ga_P:X\mapsto\Hom_{\AA}(P,X)$ induces an equivalence
of $\lan P\ran$ with the category $\mod^{fp}-R_P$ of finitely presented right modules over $R_P$.
\end{lem}

\Pf . First, 
let us construct a functor $F:\mod^{fp}-R_P\to\Qcoh(\eta_T)$. For a finitely presented module $M$
we define $F(M)$ as an object representing the functor
$$X\mapsto \Hom_{R_P}(M,\Ga_P(X)).$$ 
To see that such an object exists
we represent $M$ as the cokernel of a morphism of $R_P$-modules
$R_P^{\oplus m}\to R_P^{\oplus n}$.
Every such a morphism comes from a morphism $f:P^{\oplus m}\to P^{\oplus n}$
and one can easily see that we can take $F(M)=\coker(f)$. It is clear
from this construction that the image of
$F$ is contained in $\lan P\ran$ and that $\Ga_P(F(M))\simeq M$ for every
finitely presented $R_P$-module $M$. This implies that
$$\Hom_{R_P}(M,M')\simeq\Hom_{R_P}(M,\Ga_P(F(M')))\simeq
\Hom_{\Qcoh(\eta_T)}(F(M),F(M')),$$
so $F$ is an equivalence of $\mod^{fp}-R_P$ with the full subcategory of 
$\AA$. It is clear that the essential image of $F$ is $\lan P\ran$.
\ed

Now we can prove our main result about the category $\Qcoh^f(\eta_T)$.

\begin{thm}\label{cat-eq-thm} 
For every quasi vector bundle $P$ of finite rank the functor \eqref{functor-eq}
induces an equivalence of $\Qcoh^f(\eta_T)$ with the category of
finitely presented right modules over $R_P$.
\end{thm}

This theorem is an immediate consequence of Lemma \ref{equiv-cat-lem} and of the following
result.

\begin{prop}\label{P-cat-prop} 
For every quasi vector bundle $P$ of finite rank
the subcategory $\lan P\ran\sub\Qcoh(\eta_T)$
coincides with $\Qcoh^f(\eta_T)$.
\end{prop}

\Pf . It is clear that $\lan P\ran\sub\Qcoh^f(\eta_T)$.
Note also that $\lan P\ran$ is closed under direct sums and under passing
to direct summands. Let $Q$ be any other quasi vector bundle of finite rank.
Pick a sufficiently large number $N$ such that $N\rk P>\rk Q$ and
a quasi vector bundle $R$ of rank $N\rk P-\rk Q$. By Theorem \ref{quasi-thm}
there exists an isomorphism 
$$P^{\oplus N}\simeq Q\oplus R$$
in $\Qcoh(\eta_T)$. Hence, every quasi vector bundle of finite rank is contained in $\lan P\ran$. Since $\lan P\ran$ is closed under taking cokernels, from Theorem \ref{proj-thm} we get that
$\lan P\ran=\Qcoh^f(\eta_T)$.
\ed

\begin{rems} 1. 
We do not know whether $R_P$ is actually {\it von Neumann regular}, i.e., whether every finitely generated right ideal in it is a direct summand. 
An equivalent question is whether the category $\Qcoh^f(\eta_T)$ is semisimple.
Yet another reformulation of this question is whether every quasicoherent sheaf of finite rank is isomorphic to a quasi vector bundle in $\Qcoh(\eta_T)$ (by Theorem \ref{proj-thm}(i)).
 
\noindent
2. It is not true that $P$ is a generator of 
$\Qcoh(\eta_T)$ (even if it is a vector bundle).
More precisely, we claim that
$$\Hom_{\Qcoh(\eta_T)}(P,\oplus_{n=1}^{\infty}P)\neq\oplus_{n=1}^{\infty}
\Hom_{\Qcoh(\eta_T)}(P,P).$$
Indeed, using Lemma \ref{eps-lem} it is easy to construct a collection of nonzero subbundles
$P_n\sub P$ such that we have an embedding $\oplus_{n=1}^{\infty}P_n\sub P$
and $\sum_{n=1}^{\infty}\rk P_n=\rk P$.
Therefore, we obtain a direct sum decomposition in $\Qcoh(\eta_T)$
$$P\simeq\oplus_{n=1}^{\infty}P_n.$$
Hence, 
$$\Hom_{\Qcoh(\eta_T)}(P,\oplus_{n=1}^{\infty}P)\simeq
\Hom_{\Qcoh(\eta_T)}(\oplus_{n=1}^{\infty}P_n,\oplus_{n=1}^{\infty}P).$$
Taking an element in this space that induces an embedding of $P_n$ into the
$n$-th summand $P$, one can easily derive our claim.
\end{rems}

\end{document}